\documentclass[12pt]{article}
\usepackage{amssymb,a4}
\usepackage{amsmath,amsfonts,amssymb,amsthm}
\allowdisplaybreaks[4]
\def\deg{{\rm deg}}

\def\de{\delta}

\newcommand{\la}{\lambda}
\def\1{{\bf 1}}

\def\C{{\mathbb C}}

\def\Z{{\mathbb Z}}

\def\N{{\mathbb N}}
\def\Q{{\mathbb Q}}
\def\1{{\bf 1}}

\numberwithin{equation}{section}

\newtheorem{theorem}{Theorem}[section]
\newtheorem{proposition}[theorem]{Proposition}
\newtheorem{lemma}[theorem]{Lemma}

\newtheorem{remark}[theorem]{Remark}
\theoremstyle{definition}
\newtheorem{definition}[theorem]{Definition}

\usepackage{xparse}
\ExplSyntaxOn
\NewDocumentCommand \fun { m e{^_} }
  {%
    \operatorname{#1}%
    \IfValueT{#2}{\sp{#2}}%
    \IfValueT{#3}{\sb{#3}}%
  }
\ExplSyntaxOff

\begin{document}
\begin{center}
{\Large {\bf Bimodules associated to twisted modules of vertex operator algebras and fusion rules}} \\

\vspace{0.5cm} Yiyi Zhu
\\
Department of Mathematics\\ South China University of Technology\\ Guangzhou 510641 CHINA\\ 

\vspace{0.5cm} \today
\end{center}
\hspace{1cm}
\begin{abstract} Let $V$ be a vertex operator algebra, $T\in \N$ and $(M^k, Y_{M^k})$ for $k=1, 2, 3$ be a $g_k$-twisted module, where $g_k$ are commuting automorphisms of $V$ such that $g_k^T=1$ for $k=1, 2, 3$ and $g_3=g_1g_2$. Suppose $I(\cdot, z)$ is an intertwining operator of type $\left(\begin{array}{c}
M^{3} \\
M^{1} M^{2}
\end{array}\right)
$. We construct an $A_{g_1g_2}(V)$-$A_{g_2}(V)$-bimodule $A_{g_1g_2, g_2}(M^1)$ which determines the action of $M^1$ from the bottom level of $M^2$ to the bottom level of $M^3$ and explored its connections with fusion rules. 
\end{abstract}


\section{Introduction}
Ever since its appearance, vertex (operator) algebra has played an important role in conformal field theory \cite{MS}, and in the study of moonshine and Monster \cite{FLM}, \cite{B}. A vertex operator algebra is a vector space equipped with a linear map that sends each vector to a sequence of operators. These sequences of operators satisfy some axioms, among which the most important one is the so-called Jacobi identity. While they may look very different, vertex operator algebras can be viewed as analogs of Lie algebras and commutative associative algebras. In fact, the Jacobi identity is equivalent to associativity and commutativity \cite{FHL}. 

In this paper, we mainly focus on the associative aspect of vertex operator algebras. In \cite{Z}, the famous Zhu's algebra was constructed. Given a vertex operator algebra $V$, Zhu constructed an associative algebra $A(V)$ which is obtained from all weight-zero components of vertex operators modulo some relations hiding in Jacobi identity. Zhu in \cite{Z} also established a one-to-one correspondence between the set of equivalence classes of irreducible $A(V)$-modules and the set of equivalence classes of irreducible admissible $V$-modules. On the other hand, Zhu proved for an admissible $V$-module $M=\oplus_{n\in \N}M(n)$, the bottom level $M(0)$ is an $A(V)$-module. For an $A(V)$-module $U$, one can construct an admissible $V$-module whose bottom level is exactly $U$. In \cite{DLM98},  Zhu's construction was generalized to the twisted case. Given a vertex operator algebra $V$ and an automorphism $g$ of $V$ with finite order $T$, an associative algebra $A_g(V)$ was constructed, and the notion of $g$-twisted $V$-module was introduced. It was established in \cite{DLM98} that there is a one-to-one correspondence between the set of equivalence classes of irreducible $A_g(V)$-modules and the set of equivalence classes of irreducible admissible $g$-twisted $V$-modules. For an admissible $g$-twisted $V$-module $M=\oplus_{n\in \N}M(\frac{n}{T})$, the bottom level $M(0)$ is an $A_g(V)$-module. For an $A_g(V)$-module $U$, one can construct an admissible $g$-twisted $V$-module whose bottom level is exactly $U$. 

The Zhu's construction can also be generalized to any $V$-module, see \cite{FZ}. For this purpose, the notion of intertwining operator jumps in, see \cite{FHL}. Let $M^1$, $M^2$ and $M^3$ be three $V$-modules. An intertwining operator of type $\left(\begin{array}{c}
M^{3} \\
M^{1} M^{2}
\end{array}\right)
$ is a linear map $I$: $M^1$ $\mapsto \mbox{Hom}_{\C}(M^2, M^3)\{z\}$ satisfying similar axioms as in the definition of $V$-modules, including the Jacobi identity. Each homogeneous vector in $M^1$ corresponds to a sequence of operators in $\mbox{Hom}_\C(M^2, M^3)$. Since $\mbox{Hom}_\C(M^2(0), M^3(0))$ has an $A(V)$-$A(V)$-bimodule structure, in \cite{FZ}, Frankel and Zhu focus on weight-zero operators and constructed an $A(V)$-$A(V)$-bimodule $A(M^1)$, which is a quotient of $M^1$. As an application, they claimed that there exists a bijection between $\mbox{Hom}_{A(V)}(A(M^1)\otimes_{A(V)} M^2(0), M^3(0))$ and ${\mathcal{V}}^{M^3}_{M^1M^2}$, the space of intertwining operators of type $\left(\begin{array}{c}
M^{3} \\
M^{1} M^{2}
\end{array}\right)
$, without giving a detailed proof, see \cite{FZ}. Then in \cite{L1} and \cite{L2}, Li gave a proof when $M^2$ and $(M^3)'$ are generalized Verma $V$-module. Furthermore, a counterexample was given when $M^2$ and $(M^3)'$ are not of such type. Later, in \cite{L3}, a shorter and more conceptual proof was given using the regular representation theory developed in \cite{L4}. Theoretically, this provides us one way to compute fusion rules. In \cite{DLM4}, a series of associative algebras $A_n(V)$ for $n\geq 0$ was constructed with $A_0(V)$ being Zhu's $A(V)$. In \cite{DJ}, Dong and Jiang generalized the Frankel and Zhu's $A(V)$-$A(V)$-bimodule theory to $A_n(V)$-$A_m(V)$-bimodule theory for $m,n\geq 0$ by considering nonzero-weighted operators instead of weight-zero operators. For $m, n\geq 0$, they constructed $A_n(V)$-$A_m(V)$-bimodules $A_{n,m}(V)$ and showed how to use $A_{n,m}(V)$ to construct admissible $V$-modules from $A_m(V)$-modules. There have been many generalizations of Zhu's $A(V)$ theory (see for examples, \cite{DR}, \cite{JJ}, \cite{MT}).

The goal of this paper is to generalize the $A(V)$-$A(V)$-bimodule construction to twisted case. Let $M^i$ be a $g_i$-twisted module, $i=1, 2, 3$. There is also the notion of intertwining operator of type $\left(\begin{array}{c}
M^{3} \\
M^{1} M^{2}
\end{array}\right)
$, see \cite{X}, \cite{DLM96}. Instead of Jacobi identity, generalized Jacobi identity is required in the definition. In this paper, we shall construct an $A_{g_3}(V)$-$A_{g_2}(V)$-bimodule $A_{g_3, g_2}(M^1)$ in the case $g_3=g_1g_2$, which is necessarily true if there exists a nonzero intertwining operator of such type. In fact, given a nonzero intertwining operator $I$, one can construct an $A_{g_1g_2}(V)$-$A_{g_2}(V)$-bimodule $S_I$ and each $S_I$ is a homomorphic image of $A_{g_1g_2, g_2}(M^1)$. We try to build a bijection between ${\mathcal{V}}^{M^3}_{M^1M^2}$ and $\mbox{Hom}_{A_{g_1g_2}(V)}(A_{g_1g_2, g_2}(M^1)\otimes_{A_{g_2}(V)} M^2(0), M^3(0))$, but we are not able to establish the surjectivity in this paper. 

This paper is organized as follows: In Section 2, we recall some basics. In section 3, we first present the construction of $A_{g_1g_2, g_2}(M^1)$ and $S_I$. Then we give an $A_{g_1g_2}(V)$-$A_{g_2}(V)$-bimodule homomorphism from $A_{g_1g_2,g_2}(M^1)$ to $S_I\subset$Hom$_{\mathbb{C}}(M^2(0), M^3(0))$. Finally, we prove $$\mbox{dim}\mathcal{V}^{M^3}_{M^1M^2} \leq \mbox{dim} \mbox{Hom}_{A_{g_1g_2}(V)}(A_{g_1g_2, g_2}(M^1)\otimes_{A_{g_2}(V)} M^2(0), M^3(0)).$$

\section{Preliminary}
Throughout this paper, we denote the field of complex numbers by $\mathbb{C}$, the field of rational numbers by $\mathbb{Q}$, the ring of integers by $\mathbb{Z}$, and the set of natural numbers by $\N$. We adopt the standard notation to denote the space of $W$-valued formal series in arbitrary integer, complex power of $z$ for a vector space $W$ by $W[[z, z^{-1}]]$ and $W\{z\}$ respectively. All vector spaces are assumed to be over $\mathbb{C}$. \par
\subsection{Vertex operator algebras and modules}
Below are some definitions related to vertex (operator) algebras, see
 \cite{B}, \cite{DLM97}, \cite{FHL}, \cite{FLM}, and \cite{LL}.
 \begin{definition}
     A vertex algebra is a triple $(V, Y, \textbf{1})$, where $V$ is a vector space, $Y$: $V$ $\mapsto$ (End$V$)$[[z, z^{-1}]]$, $v\mapsto Y(v, z)=\sum\limits_{n\in \Z}v_nz^{-n-1}$, is a linear map and vacuum vector \textbf{1}$\in V$, satisfying the following conditions: 
     \begin{enumerate}
         \item[(V1).] Truncation condition: for $u, v \in V$, $Y(u, z)v=\sum\limits_{n\in \Z}u_nvz^{-n-1}$ and $u_nv=0$ for sufficiently large $n$;
    
        \item[(V2).] Vacuum property: $Y(\textbf{1}, z)=id_V$;
    
        \item[(V3).] Creation property: $Y(u, z){\1}\in V[[z]]$ and $\lim\limits_{z\to 0}Y(u, z){\1}=u$ for $u\in V$;
    
        \item[(V4).] Jacobi identity: for $u, v \in V$, $$z_0^{-1}\delta(\frac{z_1-z_2}{z_0})Y(u, z_1)Y(v, z_2)-z_0^{-1}\delta(\frac{-z_2+z_1}{z_0})Y(v, z_2)Y(u, z_1)$$
        $$=z_1^{-1}\delta(\frac{z_2+z_0}{z_1})Y(Y(u, z_0)v, z_2).$$
     \end{enumerate}
 \end{definition}
 As a consequence of Jacobi identity, we have the following associativity and commutativity. \\
 Associativity: for any $a, c \in V$, there is a nonnegative integer $r$ such that for all $b \in V$
 \[(z_0+z_2)^rY(a, z_0+z_2)Y(b, z_2)c=(z_2+z_0)^rY(Y(a, z_0)b, z_2)c.\]
 Commutativity: for any $a, b \in V$, there is a nonnegative integer $k$ such that for any $c \in V$
 \[(z_1-z_2)^kY(a, z_1)Y(b, z_2)c=(z_1-z_2)^kY(b, z_2)Y(a, z_1)c.\]
 In fact, the Jacobi identity is equivalent to the commutativity and the associativity, see \cite{FHL} and \cite{LL}.
 
\begin{definition}
    A vertex operator algebra is a quadruple $(V, Y, \textbf{1}, \omega)$, where $(V, Y, \textbf{1})$ is a vertex algebra, and $\omega$ $\in V$ satisfying the following conditions:
    \begin{enumerate}
        \item[(V5).] $V=\coprod\limits_{n\in {\Z}}V_n$, $\mbox{dim}V_n< \infty$ and $V_n=0$ if $n\ll0$;
        \item[(V6).] Virasoro relations: 
        $$[L(m), L(n)]=(m-n)L(m+n)+\frac{1}{12}(m^3-m)\delta_{m+n,0}c_V$$
        where $Y(\omega, z)=\sum_{n\in \mathbb{Z}}\omega_nz^{-n-1}=\sum_{n\in \mathbb{Z}}L(n)z^{-n-2}$, $c_V\in \mathbb{C}$;
        \item[(V7).] $L(0)$-eigenspace decomposition: $L(0)u=nu$ for $u\in V_n$. That is, the weight of $u$, denoted by $\mbox{wt} u$, is the corresponding eigenvalue of $L(0)$; 
        \item[(V8).] $L(-1)$-derivative property: $Y(L(-1)u, z)=\frac{d}{dz}Y(u, z)$;
    \end{enumerate}
\end{definition}

\begin{definition}
    An automorphism $g$ of a vertex operator algebra $V$ is a linear automorphism of $V$ such that:
    \begin{enumerate}
        \item $gY(u, z)g^{-1}=Y(gu, z)$ for all $u\in V$;
        \item $g\omega=\omega$.
    \end{enumerate}
\end{definition}

\begin{definition}
    A weak module of vertex operator algebra $V$ is a vector space $M$, equipped with a linear map $Y_M: V \mapsto$ (End$M$)$[[z^{-1}, z]]$ such that: 
    \begin{enumerate}
        \item[(M1).] $Y_M(\textbf{1}, z)=id_M$;
        \item[(M2).] For $u\in V$ and $w \in M$, $Y_M(u, z)w=\sum_{n\in \mathbb{Z}}u_nwz^{-n-1}$ and $u_nw=0$ if $n\gg0$;
        \item[(M3).] For $u, v \in V$, 
        $$z_0^{-1}\delta(\frac{z_1-z_2}{z_0})Y_M(u, z_1)Y_M(v, z_2)-z_0^{-1}\delta(\frac{-z_2+z_1}{z_0})Y_M(v, z_2)Y_M(u, z_1)$$
        $$=z_1^{-1}\delta(\frac{z_2+z_0}{z_1})Y_M(Y(u, z_0)v, z_2).$$
    \end{enumerate}
\end{definition}
\begin{definition}
    A weak $V$-module $M$ is called admissible if the following hold:
    \begin{enumerate}
        \item It is equipped with an $\N$-grading, $M=\oplus_{n\in \N} M_{n}$;
        \item For homogeneous $u\in V$, $u_nM_{m}\subset M_{m+\mbox{wt} u-n-1}$.
    \end{enumerate}
\end{definition}
That is, $M$ is an $\N$-graded space and for homogeneous $u\in V$, $u_n$ is a homogeneous map of degree $\mbox{wt} u-n-1$. In particular, $u_{\mbox{wt} u-1}$ preserves all homogeneous subspaces of $M$, which we denote by $o_M(u)$, following \cite{Z}.
 
\begin{definition}
    A weak $V$-module $M$ is called ordinary if the following hold:
    \begin{enumerate}
        \item $M=\oplus_{\la\in \C} M_{\la}$ with dim$M_{\la}<\infty$, $M_{\la+n}=0$ for all sufficiently negative integer $n$;
        \item $M_{\la}=\{w\in M \mid L(0)w=\la w\}$.
    \end{enumerate}
\end{definition}
It's straightforward to show that ordinary modules are admissible. By $V$-modules we always mean ordinary $V$-modules.  
\subsection{Twisted modules and associative algebras}
Let $V$ be a vertex operator algebra, $g$ an automorphism of $V$ with order $T<\infty$. Then $$V=\oplus_{r=0}^{T-1} V^r,$$ where $$V^r=\{v\in V\mid gv=e^{\frac{2\pi i r}{T}}v\}.$$ 
The following definitions and results can be found in \cite{DLM98}, \cite{Z}. 
\begin{definition}
    A weak $g$-twisted $V$-module $M$ is a vector space equipped with a linear map $Y_M$: $V \mapsto$ (End$M$)\{$z$\} such that: 
    \begin{enumerate}
        \item[(M1).] For $u\in V^r, w\in M$, $Y_M(u, z)w=\sum_{n\in\frac{r}{T}+\Z}u_nwz^{-n-1}$ and $u_nw=0$ if $n\gg0$;
        \item[(M2).] $Y_M(\1, z)=id_M$;
        \item[(M3).] Twisted Jacobi identity: for $u\in V^r, v \in V$,
        $$z_0^{-1}\delta(\frac{z_1-z_2}{z_0})Y_M(u, z_1)Y_M(v, z_2)-z_0^{-1}\delta(\frac{-z_2+z_1}{z_0})Y_M(v, z_2)Y_M(u, z_1)$$
        $$=z_1^{-1}\delta(\frac{z_2+z_0}{z_1})(\frac{z_2+z_0}{z_1})^{\frac{r}{T}}Y(Y(u, z_0)v, z_2).$$
    \end{enumerate}
\end{definition}
\begin{definition}
    A weak $g$-twisted $V$-module $M$ is called admissible if the following hold:
    \begin{enumerate}
        \item $M=\oplus_{n\in\frac{1}{T}\N} M_{n}$;
        \item For homogeneous $u\in V$, $u_nM_{m}\subseteq M_{m+\mbox{wt} u-n-1}$. 
    \end{enumerate}
\end{definition} 
For $u_n$ in $Y_M(u, z)$ where $n\in \frac{1}{T}\Z$, again, like in untwisted case, we set $$\mbox{wt} u_n=\mbox{wt} u-n-1$$ and set $$o_M(u)=u_{\mbox{wt} u-1}.$$ As we can see, if $M$ is admissible, then $$u_nM_{n}\subseteq M_{n+\mbox{wt} u_n}$$ and $o_M(u)$ stabilizes each homogeneous subspace $M_{n}$. 
\begin{definition}
    A $g$-twisted $V$-module $M$ is a weak $g$-twisted $V$-module carrying a $\C$-grading such that the following hold:
    \begin{enumerate}
        \item $M=\oplus_{\la\in \C} M_{\la}$ with dim$M_{\la}<\infty$, $M_{\la+\frac{n}{T}}=0$ for all sufficiently negative integers $n$;
        \item $M_{\la}=\{w\in M \mid L(0)w=\la w\}$.
    \end{enumerate}
\end{definition}
Like in untwisted case, a $g$-twisted $V$-module is admissible.\par
\begin{proposition}\label{p2.10}
Let $M$ be an irreducible $g$-twisted module, then 
\[
M=\mbox{span}\{u_na| u\in V, n\in \frac{1}{T}\mathbf{N}\}
\]
for any nonzero $a\in M$.
\end{proposition}
Following \cite{DLM98}, for $r\in \N$, define $\delta(r)=1$ if $r\equiv 0\pmod{T}$ and $\delta(r)=0$ otherwise. For $u\in V^r, v\in V$, define 
$$u\circ_gv=\mbox{Res}_z \frac{(1+z)^{\mbox{wt} u-1+\delta(r)+\frac{r}{T}}}{z^{1+\delta(r)}}Y(u, z)v.$$
Let $O_g(V)$ be the linear span of all $u\circ_gv$ and $(L(-1)+L(0))u$, define $A_g(V)$ to be the quotient space $V/O_g(V)$. 
Also define $$u\ast_gv=\begin{cases}
\mbox{Res}_z \frac{(1+z)^{\mbox{wt} u}}{z}Y_M(u,z)v, & \mbox{if $u\in V^r$ with $r=0$}\\
0, & \mbox{otherwise}
\end{cases}$$
The following results were obtained in \cite{DLM98}.
\begin{proposition}\label{p2.9}
    $V^r\subseteq O_g(V)$ for $0<r<T$.
\end{proposition}
Proposition \ref{p2.9} tells us that $A_g(V)$ is a quotient of $V^0$.
\begin{theorem}\label{t2.10}
    Let $V$ be a vertex operator algebra, $g$ an automorphism of $V$ with finite order $T$. Then $A_g(V)$ is an associative algebra with respect to the operation $\ast_g$. Furthermore, $\1+O_g(V)$ acts as the identity and $\omega+O_g(V)$ lies in the center of $A_g(V)$.
\end{theorem}
\begin{theorem}\label{t2.11}
    Let $V$ be a vertex operator algebra, $g$ an automorphism of $V$ with finite order $T$. Let $M=\oplus_{n\in\N}M_{\frac{n}{T}}$ be an admissible $g$-twisted $V$-module. Then
    \[o_M(u)o_M(v)=o_M(u\ast_g v)
    \]and 
    \[o(u')=0\]
    hold in End($M_{0}$) for every $u, v\in V$ and $u'\in O_g(V)$. Thus, the bottom level $M_{0}$ is a left $A_g(V)$-module with $u+O_g(V)$ acting as $o_{M}(u)$.
\end{theorem}
\begin{remark}
     If $g=id$, then $g$-twisted $V$-modules are just $V$-modules and $A_g(V)$ coincides with Zhu's algebra $A(V)$ which was constructed in \cite{Z}. In this special case, the untwisted version of Theorem \ref{t2.10} and \ref{t2.11} were established in \cite{Z}.
\end{remark}
\subsection{Affinization of vertex algebras and contragredient modules}
Let $V$ be a vertex operator algebra with automorphism $g$ of finite order $T$. Consider the tensor product space 
\[
\mathcal{L}(V)=\mathbb{C}[t^{\frac{1}{T}}, t^{-\frac{1}{T}}]\otimes V.
\]
Since $\mathbb{C}[t^{\frac{1}{T}}, t^{-\frac{1}{T}}]$ is a vertex algebra \cite{B}, $\mathcal{L}(V)$ is a tensor product of two vertex algebras hence a vertex algebra. $g$ can be extend to an automorphism of $\mathcal{L}(V)$ by 
\[
g(t^m\otimes u)=e^{-2\pi im}(t^m\otimes gu).
\]
Denote the $g$-invariant subspace of $\mathcal{L}(V)$ by $\mathcal{L}(V, g)$, then 
\[
\mathcal{L}(V, g)=\oplus_{r=0}^{r=T-1}t^{\frac{r}{T}}\mathbb{C}[t, t^{-1}]\otimes V^r.
\]
$\mathcal{L}(V, g)$ is a vertex subalgebra of $\mathcal{L}(V)$ with the translation operator $D=\frac{d}{dt}\otimes 1+1\otimes L(-1)$. Let 
\[
V[g]=\mathcal{L}(V, g)/D\mathcal{L}(V, g).
\]We have the following result (cf. \cite{B}, \cite{DLM98}):
\begin{proposition}
$V[g]$ has a Lie algebra structure given by 
\[
[u+D\mathcal{L}(V, g), v+D\mathcal{L}(V, g)]=u_0v+D\mathcal{L}(V, g)
\]
for $u, v \in \mathcal{L}(V, g)$.
\end{proposition}
Denote $t^m\otimes u+ D\mathcal{L}(V, g)$ by $u(n)$ for $m\in \frac{1}{T}\mathbb{Z}, u\in V$. One can introduce a $\frac{1}{T}\mathbb{Z}$-gradation on $\mathcal{L}(V)$ by defining 
\[
\mbox{deg}(t^m\otimes u)=\mbox{wt}u-m-1  
\]
for $m\in \frac{1}{T}\mathbb{Z}, u\in V$. Since $D$ is a homogeneous operator, we have 
\[
V[g]=V[g]_-\oplus V[g]_0\oplus V[g]_+,
\]
where $V[g]_{\pm}=\oplus_{0<n\in \frac{1}{T}\mathbb{Z}}V[g]_{\pm n}$.
We have the following result \cite{DLM98}:
\begin{proposition}\label{p2.16}
Let $M$ be a weak $g$-twisted module, then $u(n)\mapsto u_n$ defines a representation of Lie algebra $V[g]$ on $M$. Furthermore, if $M$ carries a $\frac{1}{T}\mathbb{N}$-grading, then $M$ is an admissible $g$-twisted module if and only if $M$ is a $\frac{1}{T}\mathbb{N}$-graded module for the graded Lie algebra $V[g]$.
\end{proposition}
Now Proposition \ref{p2.10} can be reformulated as:
\begin{proposition}\label{p2.17}
Let $M$ be an irreducible $g$-twisted module, then 
\[
M=\mbox{span}\{u(n)a| u(n)\in V[g]\}
\]
for any nonzero $a\in M$.
\end{proposition}
Let $M=\oplus_{n\in \mathbb{N}}M(\frac{n}{T})$ be a $g$-twisted module. The contragredient module $M'$ is defined as follows:
\[
M'=\oplus_{n\in \mathbb{N}}M(\frac{n}{T})^*,
\]
the vertex operator $Y_M'(u, z)$ is defined by 
\[
\langle Y_M'(u, z)f, a\rangle = \langle f, Y_{M}(e^{zL(1)}(-z^{-2})^{L(0)}u, z^{-1})a\rangle
\]
for $u\in V$, $f\in M'$ and $a\in M$. One can prove (cf. \cite{DLM98}, \cite{FHL}, \cite{X}) the following:
\begin{theorem}\label{t2.18}
    $(M', Y_M')$ is a $g^{-1}$-twisted module and $(M'', Y_M'')=(M, Y_M)$. $M$ is irreducible if and only if $M'$ is irreducible. 
\end{theorem}
Following \cite{L1} and \cite{L2}, for $u\in V^r$, set 
\[
Y_t(u, z)=\sum_{n\in \frac{r}{T}+\mathbb{Z}}u(n)z^{-n-1}\in V[g][[z^{\frac{1}{T}}, z^{-\frac{1}{T}}]].
\]
Define $\theta\in \mbox{Hom}_{\mathbb{C}}(V[g^{-1}], V[g])$ as follows:
\[
\theta(Y_t(u, z))=Y_t(e^{zL(1)}(-z^{-2})^{L(0)}u, z^{-1}).
\]
Then by Proposition \ref{p2.16},
\[
\langle u(n)f, a\rangle = \langle f, \theta(u(n))a\rangle.
\]
and by Theorem \ref{t2.18},
\[
\theta^2=1.
\]
\subsection{Intertwining operators}
From now on, we fix a branch of logarithmic function with 
\begin{equation*}
    log(z)=ln|z|+ iarg z, 0\leq arg z<2\pi, z\in \mathbb{C}^{\times}.
\end{equation*}
Then define 
$$z^\alpha=e^{log(z)\alpha}$$ for $z\in \mathbb{C}^{\times}$, $\alpha \in \mathbb{C}$.
In particular, 
$$(-1)^{\alpha}=e^{i\pi \alpha}.$$
Let $g_k$ ($k=1, 2, 3$) be three commuting automorphisms of vertex operator algebra $V$ and $T\in \N$, a finite number such that $g_k^T=1$ for $k=1, 2, 3$. Let $(M^k, Y_{M^k})$ be a $g_k$-twisted $V$-module for $k=1, 2, 3$. Since $g_1g_2=g_2g_1$, we have the following common eigenspace decomposition: 
$$V=\oplus_{0\leq j_1, j_2<T}V^{(j_1, j_2)},$$ 
where 
$$V^{(j_1, j_2)}=\{v\in V\mid g_kv=e^{\frac{2\pi ij_k}{T}}v, k=1, 2\}.$$
The following definition can be found in \cite{X}, \cite{DLM96} and \cite{LS}.
\begin{definition}
    An intertwining operator of type $\left(\begin{array}{c}
    M^{3} \\
    M^{1} M^{2}
    \end{array}\right)
    $ is a linear map $I(\cdot, z): M^1 \mapsto (\mbox{Hom}(M^2, M^3))\{z\}$ such that:
    
    (1) For $w^1\in M^1, w^2\in M^2$, $I(w^1, z)w^2=\sum_{n\in\C}w^1_nw^2z^{-n-1}$ and $w^1_{c+n}w^2=0$ for a fixed $c\in \C$, $n\gg 0$, and $n\in \Q$.  
    
    (2) Generalized Jacobi identity: for $u\in V^{(j_1,j_2)}$, $w_1 \in M^1$, $w_2\in M^2$, and $0\leq j_1, j_2\leq T-1$,
    \begin{align*}
    &z_0^{-1}\delta(\frac{z_1-z_2}{z_0})(\frac{z_1-z_2}{z_0})^{\frac{j_1}{T}}Y_{M^3}(u,z_1)I(w_1,z_2)w_2\\
    -&e^{\frac{j_1}{T}\pi i}z_0^{-1}\delta(\frac{-z_2+z_1}{z_0})(\frac{z_2-z_1}{z_0})^{\frac{j_1}{T}}I(w_1,z_2)Y_{M^2}(u,z_1)w_2\\
    =&z_1^{-1}\delta(\frac{z_2+z_0}{z_1})(\frac{z_2+z_0}{z_1})^{\frac{j_2}{T}}I(Y_{M^1}(u,z_0)w_1,z_2)w_2    
    \end{align*} 
    
    (3) $L(-1)$-derivative property: for $w_1\in M^1$, 
    \begin{align*}
    &I(L(-1)w_1, z)=\frac{d}{dz}I(w_1, z).    
    \end{align*}
\end{definition}
\begin{remark}
     Note that $Y(\cdot, z)$ acting on $V$ is an example of an intertwining operator of type $\left(\begin{array}{c}
    V \\
    V V
    \end{array}\right)
    $ and $Y_M(\cdot, z)$ acting on a $g$-twisted module $M$ is an example of an intertwining operator of type $\left(\begin{array}{c}
    M \\
    V M
    \end{array}\right)
    $.
\end{remark}
Denote the space of intertwining operators of type $\left(\begin{array}{c}
    M^{3} \\
    M^{1} M^{2}
    \end{array}\right)
$ by ${\mathcal{V}}^{M^3}_{M^1M^2}$ and set 
\begin{equation*}
N^{M^3}_{M^1M^2}=\mbox{dim}\,{\mathcal{V}}^{M^3}_{M^1M^2}.
\end{equation*}
These numbers are called fusion rules associated to these data. If \[N^{M^3}_{M^1M^2}>0,\] then \[g_3=g_1g_2\] (see \cite{X}).\par
For the rest of the paper, let $(M^k, Y_{M^k})$ ($k=1, 2, 3$) be a $g_k$-twisted $V$-module such that \[M^k=\oplus_{n\in \N} M^k_{h_k+\frac{n}{T}},\] where \[L(0)\mid_{M^k_{h_k+\frac{n}{T}}}=(h_k+\frac{n}{T})\mbox{Id},\] and  \[g_3=g_1g_2.\] 
We also require that $M^k$ is generated by $M^k_{h_k}$.
For convenience, we denote $M^k_{h_k+\frac{n}{T}}$ by $M^k(\frac{n}{T})$ and for $w\in M^k(\frac{n}{T})$, we set $$\deg w=\frac{n}{T}.$$ Let $I(\cdot, z)$ be an intertwining operator in ${\mathcal{V}}^{M^3}_{M^1M^2}$, we have the following associativity. 
\begin{proposition}\label{p2.13}
    (Associativity) For homogeneous $u\in V^{(j_1,j_2)}$, $w_1\in M^1$ and $w_2\in M^2$, we have 
    \begin{align*}
    &(z_0+z_2)^{k+\frac{j_2}{T}}Y_{M^3}(u,z_0+z_2)I(w_1,z_2)w_2\\
    =&(z_2+z_0)^{k+\frac{j_2}{T}}I(Y_{M^1}(u,z_0)w_1,z_2)w_2,
\end{align*}
where $k$ is an integer such that $z_1^{k+\frac{j_2}{T}}Y_{M^2}(u,z_1)w_2$ involves only nonnegative integer powers of $z_1$. In particular, if $w_2\in M^2(0)$, then $k$ can take value $\mbox{wt}u-1+\delta(j_2)$.
\end{proposition}
\begin{proof}
    For homogeneous $u\in V^{(j_1,j_2)}$, $w_1\in M^1$ and $w_2\in M^2(0)$, applying $\mbox{Res}_{z_1}z_1^{k+\frac{j_2}{T}}$ to the generalized Jacobi identity gives
\[
\begin{aligned}
&\mbox{Res}_{z_1}z_1^{k+\frac{j_2}{T}}z_0^{-1}\delta(\frac{z_1-z_2}{z_0})(\frac{z_1-z_2}{z_0})^{\frac{j_1}{T}}Y_{M^3}(u,z_1)I(w_1,z_2)w_2\\
=&\mbox{Res}_{z_1}z_1^{k+\frac{j_2}{T}}z_1^{-1}\delta(\frac{z_2+z_0}{z_1})(\frac{z_2+z_0}{z_1})^{\frac{j_2}{T}}I(Y_{M^1}(u,z_0)w_1,z_2)w_2.
\end{aligned}
\]
Since (see \cite{FLM})$$z_0^{-1}\delta(\frac{z_1-z_2}{z_0})(\frac{z_1-z_2}{z_0})^{\frac{j_1}{T}}=z_1^{-1}\delta(\frac{z_0+z_2}{z_1})(\frac{z_0+z_2}{z_1})^{-\frac{j_1}{T}}$$
we have
\begin{align*}
&(z_0+z_2)^{k+\frac{j_2}{T}}Y_{M^3}(u,z_0+z_2)I(w_1,z_2)w_2\\
=&(z_2+z_0)^{k+\frac{j_2}{T}}I(Y_{M^1}(u,z_0)w_1,z_2)w_2.
\end{align*}
This completes the proof.
\end{proof}\par
The following proposition is an easy corollary of the generalized Jacobi identity and $L(-1)$-derivative property. (See \cite{FHL} and \cite{FZ}.) 
\begin{proposition}
Let \[I^{\circ}(\cdot, z)=z^{h_1+h_2-h_3}I(\cdot, z).\] Then for $w_1\in M^1$, \[I^{\circ}(w_1, z)\in (\mbox{Hom}(M^2, M^3))[[z^{\frac{1}{T}}, z^{-\frac{1}{T}}]].\]
Set \[I^{\circ}(w_1, z)=\sum_{n\in \frac{1}{T}\mathbb{Z}}w_1(n)z^{-n-1},\] then for every homogeneous $w_1\in M^1$, $n\in \frac{1}{T}\mathbb{Z}$ and $m\in \frac{1}{T}\mathbb{N}$, \[w_1(n)M^2(m)\subseteq M^3(m+\deg w_1-n-1).\]
\end{proposition}
We see that $w_1(n)$ has weight $\mbox{deg}w_1-n-1$. Denote the weight zero operator $w_1(\deg w_1-1)$ by $o_I(w_1)$.
It's obvious that $I^{\circ}(\cdot, z)$ also satisfies the generalized Jacobi identity and associativity.

\section{Bimodules associated to twisted modules}
In this section, keeping the same setting as in Section 2, we will construct an $A_{g_1g_2}(V)$-$A_{g_2}(V)$-bimodule $A_{g_1g_2, g_2}(M^1)$ and explain why we define it in such way. For $r\in N$, denote the remainder of $r$ divided by $T$ by $[r]$. \par
For homogeneous $u\in V$ and $w_1\in M^1$, we define
$$u\circ_{g_1g_2,g_2} w_1=\mbox{Res}_z\frac{(1+z)^{\mbox{wt}u-1+\delta (j_2)+\frac{j_2}{T}}}{z^{1+\delta (j_1, j_2)-\frac{j_1}{T}}}Y_{M^1}(u,z)w_1,$$
where $u\in V^{(j_1, j_2)}$ and $$\delta (j_1, j_2)=
\begin{cases}
1,& \mbox{$j_2=0$}\\
1,& \mbox{$j_2\neq 0$, $j_1+j_2\geq T$}\\
0,& \mbox{$j_2\neq 0$, $j_1+j_2<T$}
\end{cases}$$
Note that $\delta (0, j_2)=\delta (j_2)$. \\
Define 
$$u*_{g_1g_2,g_2}w_1=
\begin{cases}
\mbox{Res}_z Y_{M^1}(u,z)w_1\frac{(1+z)^{\mbox{wt}u-1+\delta(j_2)+\frac{j_2}{T}}}{z^{1-\frac{j_1}{T}}}, & \mbox{$j_1+j_2\equiv 0$ (mod $T$)}\\
0, & \mbox{otherwise}
\end{cases}
$$
Define 
$$w_1*_{g_2,g_1g_2}u=
\begin{cases}
e^{-\frac{j_1}{T}\pi i}\mbox{Res}_z Y_{M^1}(u,z)w_1\frac{(1+z)^{\mbox{wt}u-1}}{z^{1-\frac{j_1}{T}}}, & \mbox{$j_2=0$}\\
0, & \mbox{otherwise} 
\end{cases}$$\\
Let $O'_{g_1g_2, g_2}(M^1)$ be the subspace of $M^1$ spanned by all $u\circ_{g_1g_2, g_2}w_1$.
\begin{remark}
    Let $g_1=g_2=1$, then $\circ_{g_1g_2, g_2}$, $\ast_{g_1g_2, g_2}$ and $\ast_{g_2, g_1g_2}$ give the same products as in \cite{FZ}, where the authors constructed an $A(V)$-$A(V)$-bimodule $A(M)$; Let $M^1=V$, then these three products give the same construction as in \cite{DLM98}, where the authors constructed an associative algebra $A_g(V)$. 
\end{remark}
\begin{lemma}\label{l3.1}
For homogeneous $u\in V$, $w_1\in M^1$, $m, n\in \Z$ and $m\geq n\geq 0$, we have
$$\mbox{Res}_z\frac{(1+z)^{\mbox{wt} u-1+\delta (j_2)+\frac{j_2}{T}+n}}{z^{1+\delta (j_1, j_2)-\frac{j_1}{T}+m}}Y_{M^1}(u,z)w_1\in O'_{g_1g_2, g_2}(M^1).$$
\end{lemma}
The proof of Lemma \ref{l3.1} is fairly standard (cf. \cite{DLM98} and \cite{Z}).
\begin{lemma}\label{l3.2}
    For $u\in V$, $u\ast_{g_1g_2, g_2}O'_{g_1g_2, g_2}(M^1)\subseteq O'_{g_1g_2,g_2}(M^1)$.
\end{lemma}
\begin{proof}
It suffices to show it holds for homogeneous $u \in V^{(j_1, j_2)}$, where $j_1+j_2\equiv 0 \pmod{T}$. Let $u \in V^{(j_1, j_2)}$, $v\in V^{(j_3, j_4)}$, $w_1\in M^1$ and $j_1+j_2\equiv 0$ (mod $T$), then 
\begin{align*}
&u\ast_{g_1g_2, g_2}\left(v\circ_{g_1g_2, g_2}w_1\right)\\
=&\mbox{Res}_{z_1}\mbox{Res}_{z_2}\frac{(1+z_1)^{\mbox{wt}u-1+\delta(j_2)+\frac{j_2}{T}}}{z_1^{1-\frac{j_1}{T}}}\frac{(1+z_2)^{\mbox{wt}v-1+\delta(j_4)+\frac{j_4}{T}}}{z_2^{1+\delta(j_3, j_4)-\frac{j_3}{T}}}Y_{M^1}(u, z_1)Y_{M^1}(v, z_2)w_1\\
=&\mbox{Res}_{z_1}\mbox{Res}_{z_2}\frac{(1+z_1)^{\mbox{wt}u-1+\delta(j_2)+\frac{j_2}{T}}}{z_1^{1-\frac{j_1}{T}}}\frac{(1+z_2)^{\mbox{wt}v-1+\delta(j_4)+\frac{j_4}{T}}}{z_2^{1+\delta(j_3, j_4)-\frac{j_3}{T}}}Y_{M^1}(v, z_2)Y_{M^1}(u, z_1)w_1\\
&+\mbox{Res}_{z_0}\mbox{Res}_{z_1}\mbox{Res}_{z_2}\frac{(1+z_1)^{\mbox{wt}u-1+\delta(j_2)+\frac{j_2}{T}}}{z_1^{1-\frac{j_1}{T}}}\frac{(1+z_2)^{\mbox{wt}v-1+\delta(j_4)+\frac{j_4}{T}}}{z_2^{1+\delta(j_3, j_4)-\frac{j_3}{T}}}z_1^{-1}\delta(\frac{z_2+z_0}{z_1})\\
&\left(\frac{z_2+z_0}{z_1}\right)^{\frac{j_1}{T}}Y_{M^1}(Y(u, z_0)v, z_2)w_1\\
\equiv &\mbox{Res}_{z_0}\mbox{Res}_{z_2}\frac{(1+z_2+z_0)^{\mbox{wt}u-1+\delta(j_2)+\frac{j_2}{T}}}{(z_2+z_0)^{1-\frac{j_1}{T}}}\frac{(1+z_2)^{\mbox{wt}v-1+\delta(j_4)+\frac{j_4}{T}}}{z_2^{1+\delta(j_3, j_4)-\frac{j_3}{T}}}\\
&Y_{M^1}(Y(u, z_0)v, z_2)w_1 \left(\mbox{mod}\ O'_{g_1g_2, g_2}(M^1)\right)\\
\equiv& \mbox{Res}_{z_0}\mbox{Res}_{z_2}\sum_{i\geq 0}{\mbox{wt}u-1+\delta(j_2)+\frac{j_2}{T} \choose i}z_0^i(1+z_2)^{\mbox{wt}u-1+\delta(j_2)+\frac{j_2}{T}-i}\sum_{j\geq 0}\\ &{-1+\frac{j_1}{T} \choose j}z_0^j z_2^{-1+\frac{j_1}{T}-j}
\frac{(1+z_2)^{\mbox{wt}v-1+\delta(j_4)+\frac{j_4}{T}}}{z_2^{1+\delta(j_3, j_4)-\frac{j_3}{T}}}Y_{M^1}(Y(u, z_0)v, z_2)w_1\\
&\left(\mbox{mod}\ O'_{g_1g_2, g_2}(M^1)\right)\\
\equiv& \mbox{Res}_{z_2}\sum_{i\geq 0}\sum_{j\geq 0}{\mbox{wt}u-1+\delta(j_2)+\frac{j_2}{T} \choose i}{-1+\frac{j_1}{T} \choose j}Y_{M^1}(u_{i+j}v, z_2)w_1 \\
& \frac{(1+z_2)^{\mbox{wt}u-1+\delta(j_2)+\frac{j_2}{T}-i+\mbox{wt}v-1+\delta(j_4)+\frac{j_4}{T}}}{z_2^{1-\frac{j_1}{T}+j+1+\delta(j_3, j_4)-\frac{j_3}{T}}}\left(\mbox{mod}\ O'_{g_1g_2, g_2}(M^1)\right)\\
\equiv& \mbox{Res}_{z_2}\sum_{i\geq 0}\sum_{j\geq 0}{\mbox{wt}u-1+\delta(j_2)+\frac{j_2}{T} \choose i}{-1+\frac{j_1}{T} \choose j}Y_{M^1}(u_{i+j}v, z_2)w_1\\ 
&\frac{(1+z_2)^{\mbox{wt}(u_{i+j}v)-1+\delta(j_2+j_4)+\frac{[j_2+j_4]}{T}+\delta(j_4)+\delta(j_2)+j-\delta(j_2+j_4)+\frac{j_2+j_4}{T}-\frac{[j_2+j_4]}{T}}}{z_2^{1+\delta([j_1+j_3], [j_2+j_4])-\frac{[j_1+j_3]}{T}+j+1+\delta(j_3, j_4)-\frac{j_1}{T}-\frac{j_3}{T}-\delta([j_1+j_3], [j_2+j_4])+\frac{[j_1+j_3]}{T}}}\\
&\left(\mbox{mod}\ O'_{g_1g_2, g_2}(M^1)\right)\\
&\equiv 0 \left(\mbox{mod}\ O'_{g_1g_2, g_2}(M^1)\right).
\end{align*}
Here we explain the last congruence. We compare 
\begin{align*}
&n=\delta(j_4)+\delta(j_2)+j-\delta(j_2+j_4)+\frac{j_2+j_4}{T}-\frac{[j_2+j_4]}{T}
\end{align*}
and 
\begin{align*}
&m=j+1+\delta(j_3, j_4)-\delta([j_1+j_3], [j_2+j_4])+\frac{[j_1+j_3]}{T}-\frac{j_1+j_3}{T}
\end{align*} case by case. Note that for $i, j\in \Z$ and $0\leq i, j<T$, $[i+j]=i+j$ when $i+j<T$ and $[i+j]=i+j-T$ when $i+j\geq T$. Since $0\leq j_k<T$ for $k=1, 2, 3, 4$ and $j_1+j_2\equiv 0\pmod{T}$, we divide it into two cases:\\
(1) If $j_1+j_2=0$, i.e. $j_1=j_2=0$, then $n=j+1$ and $m=j+1$, $m=n$;\\
(2) If $j_1+j_2=T$, then $j_1, j_2\neq0$;\par
When $j_2+j_4<T$: if $j_4=0$, then $n=j+1$ and $m=j+1$ whether $j_1+j_3<T$ or not. So $m=n$; If $j_4\neq0$, then $n=j$ and $m=j+\delta(j_3, j_4)$ when $j_1+j_3<T$ and $m=j$ when $j_1+j_3\geq T$, either way we will have $m\geq n$.\par
When $j_2+j_4=T$: then $j_4\neq 0$ and $j_3+j_4=j_1+j_3$. Hence, $\delta(j_3, j_4)-\delta([j_1+j_3], [j_2+j_4])+\frac{[j_1+j_3]}{T}-\frac{j_1+j_3}{T}=-1$. Therefore, $n=j$, $m=j$, $m=n$.\par
When $j_2+j_4>T$: then $j_4\neq 0$. Then $n=j+1$, $m=j+1$ whether $j_1+j_3<T$ or not, either case, we have $m=n$.\\
Now by Lemma \ref{l3.1}, the last congruence holds.
\end{proof}
\begin{lemma}\label{l3.3}
    For $u\in V$, $O'_{g_1g_2, g_2}(M^1)\ast_{g_2, g_1g_2}u\subseteq O'_{g_1g_2,g_2}(M^1)$.
\end{lemma}
\begin{proof}
It suffices to show it holds for homogeneous $u\in V^{(j_1, j_2)}$, where $j_2=0$. Let $u \in V^{(j_1, 0)}$, $v\in V^{(j_3, j_4)}$ and $w_1\in M^1$, then
\begin{align*}
&e^{\frac{j_1}{T}\pi i}(v\circ_{g_1g_2, g_2}w_1)*_{g_2, g_1g_2}u\\
=&\mbox{Res}_{z_1}\mbox{Res}_{z_2}\frac{(1+z_1)^{\mbox{wt}u-1}}{z_1^{1-\frac{j_1}{T}}}\frac{(1+z_2)^{\mbox{wt}v-1+\delta(j_4)+\frac{j_4}{T}}}{z_2^{1+\delta(j_3, j_4)-\frac{j_3}{T}}}Y_{M^1}(u, z_1)Y_{M^1}(v, z_2)w_1\\
=&\mbox{Res}_{z_1}\mbox{Res}_{z_2}\frac{(1+z_1)^{\mbox{wt}u-1}}{z_1^{1-\frac{j_1}{T}}}\frac{(1+z_2)^{\mbox{wt}v-1+\delta(j_4)+\frac{j_4}{T}}}{z_2^{1+\delta(j_3, j_4)-\frac{j_3}{T}}}Y_{M^1}(v, z_2)Y_{M^1}(u, z_1)w_1\\
&+\mbox{Res}_{z_0}\mbox{Res}_{z_1}\mbox{Res}_{z_2}\frac{(1+z_1)^{\mbox{wt}u-1}}{z_1^{1-\frac{j_1}{T}}}\frac{(1+z_2)^{\mbox{wt}v-1+\delta(j_4)+\frac{j_4}{T}}}{z_2^{1+\delta(j_3, j_4)-\frac{j_3}{T}}}z_1^{-1}\delta(\frac{z_2+z_0}{z_1})\\
&\left(\frac{z_2+z_0}{z_1}\right)^{\frac{j_1}{T}}Y_{M^1}(Y(u, z_0)v, z_2)w_1\\
\equiv& \mbox{Res}_{z_0}\mbox{Res}_{z_2}\frac{(1+z_2+z_0)^{\mbox{wt}u-1}}{(z_2+z_0)^{1-\frac{j_1}{T}}}\frac{(1+z_2)^{\mbox{wt}v-1+\delta(j_4)+\frac{j_4}{T}}}{z_2^{1+\delta(j_3, j_4)-\frac{j_3}{T}}}Y_{M^1}(Y(u, z_0)v, z_2)w_1 \\
&\left(\mbox{mod}\ O'_{g_1g_2, g_2}(M^1)\right)\\
\equiv& \mbox{Res}_{z_2}\sum_{i\geq 0}\sum_{j\geq 0}{\mbox{wt}u-1 \choose i}{-1+\frac{j_1}{T} \choose j}\frac{(1+z_2)^{\mbox{wt}v-1+\delta(j_4)+\frac{j_4}{T}+\mbox{wt}u-1-i}}{z_2^{1+\delta(j_3, j_4)-\frac{j_3}{T}+1-\frac{j_1}{T}+j}}\\
&Y_{M^1}(u_{i+j}v, z_2)w_1 \left(\mbox{mod}\ O'_{g_1g_2, g_2}(M^1)\right)\\
\equiv& \mbox{Res}_{z_2}\sum_{i\geq 0}\sum_{j\geq 0}{\mbox{wt}u-1 \choose i}{-1+\frac{j_1}{T} \choose j}Y_{M^1}(u_{i+j}v, z_2)w_1\\
&\frac{(1+z_2)^{\mbox{wt}(u_{i+j}v)-1+\delta(j_4)+\frac{j_4}{T}+j}}{z_2^{1+\delta([j_1+j_3], j_4)-\frac{[j_1+j_3]}{T}+\delta(j_3, j_4)-\frac{j_3}{T}+1-\frac{j_1}{T}-\delta([j_1+j_3], j_4)+\frac{[j_1+j_3]}{T}+j}}\\
&\left(\mbox{mod}\ O'_{g_1g_2, g_2}(M^1)\right)\\
\equiv& 0 \left(\mbox{mod}\ O'_{g_1g_2, g_2}(M^1)\right)
\end{align*}
Here we explain the last congruence. We compare 
\begin{align*}
&n=j
\end{align*}
and 
\begin{align*}
&m=j+\delta(j_3, j_4)+1-\delta([j_1+j_3], j_4)+\frac{[j_1+j_3]}{T}-\frac{j_1+j_3}{T}
\end{align*} case by case:\\
(1) If $j_4=0$, then $m=j+1$ when $j_1+j_3<T$ and $m=j$ when $j_1+j_3\geq T$. Either way, we have $m\geq n$;\\
(2) If $j_4\neq 0$, then $m=j+\delta(j_3, j_4)+1-\delta(j_1+j_3, j_4)\geq j$ when $j_1+j_3<T$ and $m=j+\delta(j_3, j_4)-\delta(j_1+j_3-T, j_4)\geq j$ when $j_1+j_3\geq T$. Either way, we have $m\geq n$.\\
Now by Lemma \ref{l3.1}, the last congruence holds.
\end{proof}
\begin{lemma}\label{l3.4}
For $u, v\in V$ and $w_1\in M^1$, $(u\ast_{g_1g_2, g_2}w_1)\ast_{g_2, g_1g_2}v-u\ast_{g_1g_2, g_2}(w_1\ast_{g_2, g_1g_2}v)\subseteq O'_{g_1g_2,g_2}(M^1)$.
\end{lemma}
\begin{proof}
It suffices to show it holds for homogeneous $u\in V^{(j_1, j_2)}, v\in V^{(j_3, j_4)}$, where $j_1+j_2\equiv 0 \pmod{T}$ and $j_4=0$. Let $u\in V^{(j_1, j_2)}$, $v\in V^{(j_3, 0)}$ where $j_1+j_2\equiv 0$ (mod $T$), then 
\begin{align*}
&e^{\frac{j_3}{T}\pi i}\left[(u*_{g_1g_2, g_2}w_1)*_{g_2, g_1g_2}v-u*_{g_1g_2, g_2}\left(w_1*_{g_2, g_1g_2}v\right)\right]\\
=&\mbox{Res}_{z_1}\mbox{Res}_{z_2}\frac{(1+z_1)^{\mbox{wt}u-1+\delta(j_2)+\frac{j_2}{T}}}{z_1^{1-\frac{j_1}{T}}}\frac{(1+z_2)^{\mbox{wt}v-1}}{z_2^{1-\frac{j_3}{T}}}Y_{M^1}(v, z_2)Y_{M^1}(u, z_1)w_1\\
&-\mbox{Res}_{z_1}\mbox{Res}_{z_2}\frac{(1+z_1)^{\mbox{wt}u-1+\delta(j_2)+\frac{j_2}{T}}}{z_1^{1-\frac{j_1}{T}}}\frac{(1+z_2)^{\mbox{wt}v-1}}{z_2^{1-\frac{j_3}{T}}}Y_{M^1}(u, z_1)Y_{M^1}(v, z_2)w_1\\
=&-\mbox{Res}_{z_0}\mbox{Res}_{z_1}\mbox{Res}_{z_2}\frac{(1+z_1)^{\mbox{wt}u-1+\delta(j_2)+\frac{j_2}{T}}}{z_1^{1-\frac{j_1}{T}}}\frac{(1+z_2)^{\mbox{wt}v-1}}{z_2^{1-\frac{j_3}{T}}}z_1^{-1}\delta(\frac{z_2+z_0}{z_1})\\
&\left(\frac{z_2+z_0}{z_1}\right)^{\frac{j_1}{T}}Y_{M^1}(Y(u, z_0)v, z_2)w_1\\
=&-\mbox{Res}_{z_0}\mbox{Res}_{z_1}\mbox{Res}_{z_2}\frac{(1+z_2+z_0)^{\mbox{wt}u-1+\delta(j_2)+\frac{j_2}{T}}}{(z_2+z_0)^{1-\frac{j_1}{T}}}\frac{(1+z_2)^{\mbox{wt}v-1}}{z_2^{1-\frac{j_3}{T}}}z_1^{-1}\delta(\frac{z_2+z_0}{z_1})\\
&Y_{M^1}(Y(u, z_0)v, z_2)w_1\\
=&-\mbox{Res}_{z_0}\mbox{Res}_{z_2}\frac{(1+z_2+z_0)^{\mbox{wt}u-1+\delta(j_2)+\frac{j_2}{T}}}{(z_2+z_0)^{1-\frac{j_1}{T}}}\frac{(1+z_2)^{\mbox{wt}v-1}}{z_2^{1-\frac{j_3}{T}}}Y_{M^1}(Y(u, z_0)v, z_2)w_1\\
=&-\mbox{Res}_{z_0}\mbox{Res}_{z_2}\sum_{i\geq 0}{\mbox{wt}u-1+\delta(j_2)+\frac{j_2}{T} \choose i}z_0^{i}(1+z_2)^{\mbox{wt}u-1+\delta(j_2)+\frac{j_2}{T}-i}\\
&\sum_{j\geq 0}{-1+\frac{j_1}{T} \choose j}z_0^j z_2^{-1+\frac{j_1}{T}-j}
\frac{(1+z_2)^{\mbox{wt}v-1}}{z_2^{1-\frac{j_3}{T}}}Y_{M^1}(Y(u, z_0)v, z_2)w_1\\
=&-\mbox{Res}_{z_2}\sum_{i\geq 0}\sum_{j\geq 0}{\mbox{wt}u-1+\delta(j_2)+\frac{j_2}{T} \choose i}{-1+\frac{j_1}{T} \choose j}\frac{(1+z_2)^{\mbox{wt}v-1+\mbox{wt}u-1+\delta(j_2)+\frac{j_2}{T}-i}}{z_2^{2-\frac{j_1}{T}-\frac{j_3}{T}+j}}\\
&Y_{M^1}(u_{i+j}v, z_2)w_1\\
=&-\mbox{Res}_{z_2}\sum_{i\geq 0}\sum_{j\geq 0}{\mbox{wt}u-1+\delta(j_2)+\frac{j_2}{T} \choose i}{-1+\frac{j_1}{T} \choose j}\frac{(1+z_2)^{\mbox{wt}(u_{i+j}v)-1+\delta(j_2)+\frac{j_2}{T}+j}}{z_2^{2-\frac{j_1+j_3}{T}+j}}\\
&Y_{M^1}(u_{i+j}v, z_2)w_1\\
&\in O'_{g_1g_2, g_2}(M^1).
\end{align*} 
Again, we used Lemma \ref{l3.1} in the last step.
\end{proof}\par
To construct the desired bimodule, we need to mod out a bigger subspace than $O'_{g_1g_2, g_2}(M^1)$ from $M^1$. Let $O''_{g_1g_2, g_2}(M^1)$ be the linear span of all $(u\ast_{g_1g_2}v)\ast_{g_1g_2, g_2}w_1-u\ast_{g_1g_2, g_2}(v\ast_{g_1g_2, g_2}w_1)$, $w_1\ast_{g_2, g_1g_2}(v\ast_{g_2}u)-(w_1\ast_{g_2, g_1g_2}v)\ast_{g_2, g_1g_2}u$, $u'\ast_{g_1g_2, g_2}w_1$ and $w_1\ast_{g_2, g_1g_2}v'$, where $u, v \in V$, $u'\in O_{g_1g_2}(V)$, $v'\in O_{g_2}(V)$ and $w_1\in M^1$. Let 
$$O_{g_1g_2, g_2}(M^1)=O'_{g_1g_2, g_2}(M^1)+O''_{g_1g_2, g_2}(M^1),$$ and $$A_{g_1g_2, g_2}(M^1)=M^1/O_{g_1g_2, g_2}(M^1).$$
\begin{lemma}\label{l3.5}
For $a\in V$, $a\ast_{g_1g_2, g_2}O_{g_1g_2, g_2}(M^1)\subseteq O_{g_1g_2, g_2}(M^1)$ and $O_{g_1g_2, g_2}(M^1)\ast_{g_2, g_1g_2}a\subseteq O_{g_1g_2, g_2}(M^1)$.
\end{lemma}
\begin{proof}
It suffices to show $a\ast_{g_1g_2, g_2}O''_{g_1g_2, g_2}(M^1)\subseteq O_{g_1g_2, g_2}(M^1)$ and $O''_{g_1g_2, g_2}(M^1)\ast_{g_2, g_1g_2}a\subseteq O_{g_1g_2, g_2}(M^1)$ due to Lemmas \ref{l3.2} and \ref{l3.3}. We verify it for all 4 types of spanning vectors in $O''_{g_1g_2, g_2}(M^1)$. \\
For $u, v\in V$ and $w_1\in M^1$,
\begin{align*}
    &a\ast_{g_1g_2, g_2}\big((u\ast_{g_1g_2}v)\ast_{g_1g_2, g_2}w_1-u\ast_{g_1g_2, g_2}(v\ast_{g_1g_2, g_2}w_1)\big)\\
    =&a\ast_{g_1g_2, g_2}\big((u\ast_{g_1g_2}v)\ast_{g_1g_2, g_2}w_1\big)-a\ast_{g_1g_2, g_2}\big(u\ast_{g_1g_2, g_2}(v\ast_{g_1g_2, g_2}w_1)\big)\\
    \intertext{by the definition of $O''_{g_1g_2, g_2}(M^1)$,}
    \in&\big(a\ast_{g_1g_2}(u\ast_{g_1g_2}v)\big)\ast_{g_1g_2, g_2}w_1-\big(a\ast_{g_1g_2}u\big)\ast_{g_1g_2, g_2}(v\ast_{g_1g_2, g_2}w_1)+O''_{g_1g_2, g_2}(M^1)\\
    =&\big((a\ast_{g_1g_2}u)\ast_{g_1g_2}v\big)\ast_{g_1g_2, g_2}w_1-\big(a\ast_{g_1g_2}u\big)\ast_{g_1g_2, g_2}(v\ast_{g_1g_2, g_2}w_1)+O''_{g_1g_2, g_2}(M^1)\\
    =&\big(a\ast_{g_1g_2}u\big)\ast_{g_1g_2, g_2}(v\ast_{g_1g_2, g_2}w_1)-\big(a\ast_{g_1g_2}u\big)\ast_{g_1g_2, g_2}(v\ast_{g_1g_2, g_2}w_1)+O''_{g_1g_2, g_2}(M^1)\\
    =&O''_{g_1g_2, g_2}(M^1)\subseteq O_{g_1g_2, g_2}(M^1).
\end{align*}
The second equality holds because $A_{g_1g_2}(V)$ is associative, see Theorem \ref{t2.10}, thus
\[a\ast_{g_1g_2}(u\ast_{g_1g_2}v)\in (a\ast_{g_1g_2}u)\ast_{g_1g_2}v+O_{g_1g_2}(V),\] and \[O_{g_1g_2}(V)\ast_{g_1g_2, g_2}w_1\in O''_{g_1g_2, g_2}(M^1)\] by the definition of $O''_{g_1g_2, g_2}(M^1)$.

For $u, v\in V$ and $w_1\in M^1$,
\begin{align*}
    &\big((u\ast_{g_1g_2}v)\ast_{g_1g_2, g_2}w_1-u\ast_{g_1g_2, g_2}(v\ast_{g_1g_2, g_2}w_1)\big)\ast_{g_2, g_1g_2}a\\
    =&\big((u\ast_{g_1g_2}v)\ast_{g_1g_2, g_2}w_1\big)\ast_{g_2,g_1g_2}a-\big(u\ast_{g_1g_2, g_2}(v\ast_{g_1g_2, g_2}w_1)\big)\ast_{g_2, g_1g_2}a
    &\intertext{applying Lemma \ref{l3.4} to both terms,} 
    \in&(u\ast_{g_1g_2}v)\ast_{g_1g_2, g_2}\big(w_1\ast_{g_2, g_1g_2}a\big)-u\ast_{g_1g_2, g_2}\big((v\ast_{g_1g_2, g_2}w_1)\ast_{g_2, g_1g_2}a\big)+O'_{g_1g_2, g_2}(M^1)
    \intertext{by the definition of $O''_{g_1g_2, g_2}(M^1)$ and applying Lemma \ref{l3.4} to the second term,} 
    \in&u\ast_{g_1g_2, g_2}\big(v\ast_{g_1g_2, g_2}(w_1\ast_{g_2, g_1g_2}a)\big)-u\ast_{g_1g_2, g_2}\big(v\ast_{g_1g_2, g_2}(w_1\ast_{g_2, g_1g_2}a)\big)\\
    &+u\ast_{g_1g_2, g_2}O'_{g_1g_2, g_2}(M^1)+O'_{g_1g_2, g_2}(M^1)+O''_{g_1g_2, g_2}(M^1)\\
    \intertext{by Lemma \ref {l3.2},}
    \subseteq&O_{g_1g_2, g_2}(M^1).
\end{align*}
Similarly, one can prove 
\[a\ast_{g_1g_2, g_2}\big(w_1\ast_{g_2, g_1g_2}(v\ast_{g_2}u)-(w_1\ast_{g_2, g_1g_2}v)\ast_{g_2, g_1g_2}u\big)\in O_{g_1g_2, g_2}(M^1)\] and 
\[\big(w_1\ast_{g_2, g_1g_2}(v\ast_{g_2}u)-(w_1\ast_{g_2, g_1g_2}v)\ast_{g_2, g_1g_2}u\big)\ast_{g_2, g_1g_2}a\in O_{g_1g_2, g_2}(M^1).\]
For $u'\in O_{g_1g_2}(V)$ and $w_1\in M^1$,
\begin{align*}
    &a\ast_{g_1g_2, g_2}(u'\ast_{g_1g_2, g_2}w_1)\\
    \in&(a\ast_{g_1g_2}u')\ast_{g_1g_2, g_2}w_1+O''_{g_1g_2, g_2}(M^1)\\
    \intertext{since $O_{g_1g_2}(V)$ is an ideal of V with respect to $\ast_{g_1g_2}$,}
    \subseteq&O_{g_1g_2}(V)\ast_{g_1g_2, g_2}w_1+O''_{g_1g_2, g_2}(M^1)\\
    \subseteq&O''_{g_1g_2, g_2}(M^1)\subseteq O_{g_1g_2, g_2}(M^1).
\end{align*}
For $v'\in O_{g_2}(V)$ and $w_1\in M^1$,
\begin{align*}
    &(u'\ast_{g_1g_2, g_2}w_1)\ast_{g_2, g_1g_2}a
    \intertext{by Lemma \ref{l3.4},}
    \in&u'\ast_{g_1g_2, g_2}(w_1\ast_{g_2, g_1g_2}a)+O'_{g_1g_2, g_2}(M^1)\\
    \subseteq&O''_{g_1g_2, g_2}(M^1)+O'_{g_1g_2, g_2}(M^1)\\
    =&O_{g_1g_2, g_2}(M^1)
\end{align*}
Similarly, we can prove 
$$a\ast_{g_1g_2, g_2}(w_1\ast_{g_2, g_1g_2}v')\in O_{g_1g_2, g_2}(M^1)$$
and $$(w_1\ast_{g_2, g_1g_2}v')\ast_{g_2, g_1g_2}a\in O_{g_1g_2, g_2}(M^1).$$
Now, the proof is completed. 
\end{proof}
\begin{theorem}
$A_{g_1g_2, g_2}(M^1)$ is an $A_{g_1g_2}(V)$-$A_{g_2}(V)$-bimodule with left action $\ast_{g_1g_2, g_2}$ and right action $\ast_{g_2, g_1g_2}$.  
\end{theorem}
\begin{proof}
Combining Lemmas \ref{l3.2}, \ref{l3.3}, \ref{l3.4} and \ref{l3.5}, we see it immediately.
\end{proof}
\begin{remark}\label{r3.7}
     Consider two special cases. (1): $g_1=g_2=1$. This special case was dealt with in \cite{FZ}. The $A(V)$-$A(V)$-bimodule $A(M^1)$ constructed in \cite{FZ} is just $A_{1, 1}(M^1)$; (2): $M^1=V$. In this case, $g_1=1$ and any intertwining operator $I(\cdot, z)\in {\mathcal{V}}^{M^3}_{M^1M^2}$ is just a $g_2$-twisted $V$-module map. The associative algebra $A_{g_2}(V)$ constructed in \cite{DLM98} is just $A_{g_2, g_2}(V)$. It's worth to point out that in these two special cases, they were both able to prove that $O_{g_1g_2, g_2}(M^1)=O'_{g_1g_2, g_2}(M^1)$. It's reasonable to conjecture that $O_{g_1g_2, g_2}(M^1)=O'_{g_1g_2, g_2}(M^1)$ holds in general, but we are not able to prove it in this paper. 
\end{remark}
Next, we are going to explain our bimodule construction by connecting it to representation theory. Suppose $I(\cdot, z) \in {\mathcal{V}}^{M^3}_{M^1M^2}$. Recall that $I^{\circ}(w_1, z)=\sum_{n\in \frac{1}{T}\Z}w_1(n)z^{-n-1}$ and $o_I(w_1)=w_1(\deg w_1-1)$.
\begin{lemma}\label{l3.8}
For homogeneous $u\in V^{(j_1, j_2)}$, $w_1\in M^1$ and $w_2\in M^2(0)$, $o_I(u\ast_{g_1g_2, g_2}w_1)w_2=o_{M^3}(u)o_I(w_1)w_2$, $o_I(w_1\ast_{g_2, g_1g_2}u)w_2=o_I(w_1)o_{M^2}(u)w_2$. 
\end{lemma}
\begin{proof}
By the definition of twisted module, action $\ast_{g_1g_2, g_2}$, and action $\ast_{g_2, g_1g_2}$:\\
if $$j_1+j_2\not\equiv 0\pmod{T},$$ then $$o_I(u\ast_{g_1g_2, g_2}w_1)w_2=0=o_{M^3}(u)o_I(w_1)w_2;$$ If $$j_2\neq 0,$$ then $$o_I(w_1\ast_{g_2, g_1g_2}u)w_2=0=o_I(w_1)o_{M^2}(u)w_2.$$
So it suffices to prove the two identities for $j_1+j_2\equiv 0 \pmod{T}$ and $j_2=0$.\par
With the help of associativity, for homogeneous $u\in V^{(j_1, j_2)}$, $j_1+j_2\equiv 0\pmod{T}$, homogeneous $w_1\in M^1$ and $w_2\in M^2(0)$, we have 
\begin{align*}
&o_{M^3}(u)o_I(w_1)w_2\\
=&\mbox{Res}_{z_1}\mbox{Res}_{z_2}z_1^{\mbox{wt}u-1}z_2^{\mbox{deg}w_1-1}Y_{M^3}(u,z_1)I^{\circ}(w_1,z_2)w_2\\
=&\mbox{Res}_{z_0}\mbox{Res}_{z_1}\mbox{Res}_{z_2}z_0^{-1}\delta(\frac{z_1-z_2}{z_0})z_1^{\mbox{wt}u-1}z_2^{\mbox{deg}w_1-1}Y_{M^3}(u,z_1)I^{\circ}(w_1,z_2)w_2\\
=&\mbox{Res}_{z_0}\mbox{Res}_{z_1}\mbox{Res}_{z_2}z_1^{-1}\delta(\frac{z_0+z_2}{z_1})z_1^{\mbox{wt}u-1}z_2^{\mbox{deg}w_1-1}Y_{M^3}(u,z_1)I^{\circ}(w_1,z_2)w_2\\
=&\mbox{Res}_{z_0}\mbox{Res}_{z_2}(z_0+z_2)^{\mbox{wt}u-1}z_2^{\mbox{deg}w_1-1}Y_{M^3}(u,z_0+z_2)I^{\circ}(w_1,z_2)w_2\\
=&\mbox{Res}_{z_0}\mbox{Res}_{z_2}z_2^{\mbox{deg}w_1-1}(z_0+z_2)^{-\delta(j_2)-\frac{j_2}{T}}(z_0+z_2)^{\mbox{wt}u-1+\delta(j_2)+\frac{j_2}{T}}Y_{M^3}(u,z_0+z_2)I^{\circ}(w_1,z_2)w_2\\
&\mbox{Expanding}\ (z_0+z_2)^{-\delta(j_2)-\frac{j_2}{T}}, \mbox{we can see that only the first term}\ z_0^{-\delta(j_2)-\frac{j_2}{T}} \mbox{remains}\\ &\mbox{after applying Res}_{z_2}.\\
=&\mbox{Res}_{z_0}\mbox{Res}_{z_2}z_2^{\mbox{deg}w_1-1}z_0^{-\delta(j_2)-\frac{j_2}{T}}(z_0+z_2)^{\mbox{wt}u-1+\delta(j_2)+\frac{j_2}{T}}Y_{M^3}(u,z_0+z_2)I^{\circ}(w_1,z_2)w_2\\
=&\mbox{Res}_{z_0}\mbox{Res}_{z_2}z_2^{\mbox{deg}w_1-1}z_0^{-\delta(j_2)-\frac{j_2}{T}}(z_2+z_0)^{\mbox{wt}u-1+\delta(j_2)+\frac{j_2}{T}}I^{\circ}(Y_{M^1}(u,z_0)w_1,z_2)w_2\\
=&\mbox{Res}_{z_0}\mbox{Res}_{z_2}\sum_{i\geq 0}{\mbox{wt}u-1+\delta(j_2)+\frac{j_2}{T}\choose i}z_2^{\mbox{deg}w_1-1+\mbox{wt}u-1+\delta(j_2)+\frac{j_2}{T}-i}z_0^{-\delta(j_2)-\frac{j_2}{T}+i}\\
&I^{\circ}(Y_{M^1}(u,z_0)w_1,z_2)w_2\\
=&\mbox{Res}_{z_2}\sum_{i\geq 0}{\mbox{wt}u-1+\delta(j_2)+\frac{j_2}{T}\choose i}z_2^{\mbox{deg}w_1-1+\mbox{wt}u-1+\delta(j_2)+\frac{j_2}{T}-i}I^{\circ}(u_{i-\delta(j_2)-\frac{j_2}{T}}w_1,z_2)w_2\\
=&o_I\left(\sum_{i\geq 0}{\mbox{wt}u-1+\delta(j_2)+\frac{j_2}{T}\choose i}u_{i-\delta(j_2)-\frac{j_2}{T}}w_1\right)w_2\\
=&o_I\left(\mbox{Res}_z Y_{M^1}(u,z)w_1\frac{(1+z)^{\mbox{wt}u-1+\delta(j_2)+\frac{j_2}{T}}}{z^{\delta(j_2)+\frac{j_2}{T}}}\right)w_2\\
=&o_I\left(\mbox{Res}_z Y_{M^1}(u,z)w_1\frac{(1+z)^{\mbox{wt}u-1+\delta(j_2)+\frac{j_2}{T}}}{z^{1-\frac{j_1}{T}}}\right)w_2 
\end{align*}
The last equality holds because when $j_1+j_2\equiv0$ (mod T), either $j_1=j_2=0$ or $j_1+j_2=T$, so $\delta(j_2)+\frac{j_2}{T}=1-\frac{j_1}{T}$.\par
For homogeneous $u\in V^{(j_1, 0)}$, homogeneous $w_1\in M^1$ and $w_2\in M^2(0)$, we have 
\begin{align*}
&e^{\frac{j_1}{T}\pi i}o_I(w_1*_{g_2,g_1g_2}u)w_2\\
=&o_I(\sum_{i\geq 0}{\mbox{wt}u-1\choose i}u_{i-1+\frac{j_1}{T}}w_1)w_2\\
=&\sum_{i\geq 0}{\mbox{wt}u-1\choose i}(u_{i-1+\frac{j_1}{T}}w_1)(\mbox{deg}w_1+\mbox{wt}u-i-\frac{j_1}{T}-1)w_2\\
=&\mbox{Res}_{z_2}\sum_{i\geq 0}{\mbox{wt}u-1\choose i}z_2^{\mbox{deg}w_1+\mbox{wt}u-i-\frac{j_1}{T}-1}I^{\circ}(u_{i-1+\frac{j_1}{T}}w_1,z_2)w_2\\
=&\mbox{Res}_{z_0}\mbox{Res}_{z_2}\sum_{i\geq 0}{\mbox{wt}u-1\choose i}z_2^{\mbox{deg}w_1+\mbox{wt}u-i-\frac{j_1}{T}-1}I^{\circ}(z_0^{i-1+\frac{j_1}{T}}Y_{M^1}(u,z_0)w_1,z_2)w_2\\
=&\mbox{Res}_{z_0}\mbox{Res}_{z_2}(z_2+z_0)^{\mbox{wt}u-1}z_2^{\mbox{deg}w_1-\frac{j_1}{T}}z_0^{-1+\frac{j_1}{T}}I^{\circ}(Y_{M^1}(u,z_0)w_1,z_2)w_2\\
=&\mbox{Res}_{z_0}\mbox{Res}_{z_1}\mbox{Res}_{z_2}(z_2+z_0)^{\mbox{wt}u-1}z_2^{\mbox{deg}w_1-\frac{j_1}{T}}z_0^{-1+\frac{j_1}{T}}z_1^{-1}\delta(\frac{z_2+z_0}{z_1})I^{\circ}(Y_{M^1}(u,z_0)w_1,z_2)w_2\\
=&\mbox{Res}_{z_0}\mbox{Res}_{z_1}\mbox{Res}_{z_2}z_1^{\mbox{wt}u-1}z_2^{\mbox{deg}w_1-\frac{j_1}{T}}z_0^{-1+\frac{j_1}{T}}z_1^{-1}\delta(\frac{z_2+z_0}{z_1})I^{\circ}(Y_{M^1}(u,z_0)w_1,z_2)w_2\\
=&\mbox{Res}_{z_0}\mbox{Res}_{z_1}\mbox{Res}_{z_2}z_1^{\mbox{wt}u-1}z_2^{\mbox{deg}w_1-\frac{j_1}{T}}z_0^{-1+\frac{j_1}{T}}z_0^{-1}\delta(\frac{z_1-z_2}{z_0})(\frac{z_1-z_2}{z_0})^{\frac{j_1}{T}}Y_{M^3}(u,z_1)I^{\circ}(w_1,z_2)w_2\\
&-\mbox{Res}_{z_0}\mbox{Res}_{z_1}\mbox{Res}_{z_2}e^{\frac{j_1}{T}\pi i}z_1^{\mbox{wt}u-1}z_2^{\mbox{deg}w_1-\frac{j_1}{T}}z_0^{-1+\frac{j_1}{T}}z_0^{-1}\delta(\frac{-z_2+z_1}{z_0})(\frac{z_2-z_1}{z_0})^\frac{j_1}{T}I^{\circ}(w_1,z_2)Y_{M^2}(u,z_1)w_2\\
=&\mbox{Res}_{z_1}\mbox{Res}_{z_2}e^{\frac{j_1}{T}\pi i}z_1^{\mbox{wt}u-1}z_2^{\mbox{deg}w_1-\frac{j_1}{T}}(z_2-z_1)^{-1+\frac{j_1}{T}}I^{\circ}(w_1,z_2)Y_{M^2}(u,z_1)w_2\\
=&e^{\frac{j_1}{T}\pi i}o_I(w_1)o_{M^2}(u)w_2.\\
\end{align*}
The last equality holds because $z_1^{\mbox{wt} u}Y_{M^2}(u,z_1)w_2\in M^2[[z_1^{\frac{1}{T}}]]$.    
\end{proof}
\begin{proposition}\label{l3.9}
For all $w \in O_{g_1g_2, g_2}(M^1)$, $o_I(w)|_{M^2(0)}=0$.
\end{proposition} 
\begin{proof}
First we prove it for $w\in O'_{g_1g_2, g_2}(M^1)$. Let $u\in V^{(j_1,j_2)}$ be homogeneous, $w_1\in M^1$ and $w_2\in M^2(0)$. 
\begin{align*}
    &o_I(u\circ_{g_1g_2, g_2}w_1)w_2\\
    =&o_I\big(\mbox{Res}_z\frac{(1+z)^{\mbox{wt}u-1+\delta (j_2)+\frac{j_2}{T}}}{z^{1+\delta (j_1, j_2)-\frac{j_1}{T}}}Y_{M^1}(u,z)w_1\big)w_2\\
    =&\sum_{i\geq 0}{\mbox{wt} u-1+\de(j_2)+\frac{j_2}{T}\choose i}\\
    &(u_{i-1-\delta(j_1, j_2)+\frac{j_1}{T}}w_1)(\deg w_1+\mbox{wt} u-i+\delta(j_1, j_2)-\frac{j_1}{T}-1)w_2\\
    =&\mbox{Res}_{z_2}\sum_{i\geq 0}{\mbox{wt} u-1+\delta(j_2)+\frac{j_2}{T}\choose i}z_2^{\deg w_1+\mbox{wt} u-i+\delta(j_1, j_2)-\frac{j_1}{T}-1}\\
    &I^{\circ}(u_{i-1-\delta(j_1, j_2)+\frac{j_1}{T}}w_1,z_2)w_2\\
    =&\mbox{Res}_{z_0}\mbox{Res}_{z_2}\sum_{i\geq 0}{\mbox{wt} u-1+\delta(j_2)+\frac{j_2}{T}\choose i}z_0^{i-1-\delta(j_1, j_2)+\frac{j_1}{T}}z_2^{\deg w_1+\mbox{wt} u-i+\delta(j_1, j_2)-\frac{j_1}{T}-1}\\
    &I^{\circ}(Y_{M^1}(u,z_0)w_1,z_2)w_2\\
    =&\mbox{Res}_{z_0}\mbox{Res}_{z_2}z_0^{-1-\delta(j_1, j_2)+\frac{j_1}{T}}z_2^{\deg w_1+\delta(j_1, j_2)-\delta(j_2)-\frac{j_1+j_2}{T}}(z_2+z_0)^{\mbox{wt}u-1+\delta(j_2)+\frac{j_2}{T}}\\
    &I^{\circ}(Y_{M^1}(u,z_0)w_1,z_2)w_2\\
    =&\mbox{Res}_{z_0}\mbox{Res}_{z_2}z_0^{-1-\delta(j_1, j_2)+\frac{j_1}{T}}z_2^{\deg w_1+\delta(j_1, j_2)-\delta(j_2)-\frac{j_1+j_2}{T}}(z_0+z_2)^{\mbox{wt} u-1+\delta(j_2)+\frac{j_2}{T}}\\
    &Y_{M^3}(u,z_0+z_2)I^{\circ}(w_1,z_2)w_2\\
\end{align*}
Since $z_2^{\mbox{deg} w_1}I^{\circ}(w_1,z_2)w_2\in M^3[[z_2^{\frac{1}{T}}]]$ and $\delta(j_1, j_2)-\delta(j_2)-\frac{j_1+j_2}{T}>-1$, the last expression becomes 0 after applying $\mbox{Res}_{z_2}$. \\
For those vectors in $O''_{g_1g_2, g_2}(M^1)$, we prove it for one case:
\[o_I\big((u\ast_{g_1g_2}v)\ast_{g_1g_2, g_2}w_1-u\ast_{g_1g_2, g_2}(v\ast_{g_1g_2, g_2}w_1)\big)\mid_{M^2(0)}=0,\]
the proof for other cases are similar. Use Lemma \ref{l3.8} and Remark \ref{r3.7},
\begin{align*}
    &o_I\big((u\ast_{g_1g_2}v)\ast_{g_1g_2, g_2}w_1-u\ast_{g_1g_2, g_2}(v\ast_{g_1g_2, g_2}w_1\big))\mid_{M^2(0)}\\
    =&o_{M^3}(u\ast_{g_1g_2}v)o_I(w_1)\mid_{M^2(0)}-o_{M^3}(u)o_I(v\ast_{g_1g_2, g_2}w_1)\mid_{M^2(0)}\\
    =&o_{M^3}(u)o_{M^3}(v)o_I(w_1)\mid_{M^2(0)}-o_{M^3}(u)o_{M^3}(v)o_I(w_1)\mid_{M^2(0)}\\
    =&0
\end{align*}
This completes the proof.
\end{proof}

By Theorem \ref{t2.11}, $M^2(0)$ is a left $A_{g_2}(V)$ module and $M^3(0)$ is a left $A_{g_1g_2}(V)$ module, hence $\mbox{Hom}_\mathbb{C}(M^2(0), M^3(0))$ would be an $A_{g_1g_2}(V)-A_{g_2}(V)$ bimodule with the following left and right actions:
$$\big(u+O_{g_1g_2}(V)\big)\cdot f=o_{M^3}(u)\circ f, $$ 
$$f\cdot \big(v+O_{g_2}(V)\big)=f\circ o_{M^2}(v)$$ where $f\in \mbox{Hom}_\mathbb{C}(M^2(0), M^3(0))$, $u+O_{g_1g_2}(V)\in A_{g_1g_2}(V)$ and $v+O_{g_2}(V)\in A_{g_2}(V)$. Consider the set
$$S_I:=\{o_I(w_1)|_{M^2(0)}: w_1\in M^1\}.$$ It is a subspace of $\mbox{Hom}_\mathbb{C}(M^2(0), M^3(0))$. Lemma \ref{l3.8} tells us that $S_I$ is actually a subbimodule of $\mbox{Hom}_\C(M^2(0), M^3(0))$. Regarding $o_I$ as a linear map from $M^1$
to $S_I$, we obtain that $M^1/ker$ $o_I$ $\cong S_I$ also has an $A_{g_1g_2}(V)$-$A_{g_2}(V)$-bimodule structure. From Lemma \ref{l3.9}, we see that 
\begin{proposition}\label{p3.11}
For every intertwining operator $I(\cdot, z) \in {\mathcal{V}}^{M^3}_{M^1M^2}$, there exists an $A_{g_1g_2}(V)$-$A_{g_2}(V)$-bimodule epimorphism from $A_{g_1g_2, g_2}(M^1)$ to $S_I$.
\end{proposition}
\begin{proof}
By Proposition \ref{l3.9}, $O_{g_1g_2, g_2}(M^1)\subseteq ker$ $o_I$. The statement follows immediately.
\end{proof}
\begin{remark}
Though not a perfect explanation, Proposition \ref{l3.9} and \ref{p3.11} do give us a clue why we should mod $O_{g_1g_2, g_2}(M^1)$ out. We have a series of $A_{g_1g_2}(V)-A_{g_2}(V)-$bimodules, i.e. $S_I$, where $I$ ranges through all intertwining operators in ${\mathcal{V}}^{M^3}_{M^1M^2}$. But these $S_I$'s are not good enough, because they rely on the choice of $I$. We want something that is universal or at least independent of the choice of $I$. Proposition \ref{p3.11} makes $A_{g_1g_2, g_2}(M^1)$ a good candidate.   
\end{remark}

Now we present a result about fusion rules. \par

For an intertwining operator $I\in \mathcal{V}^{M^3}_{M^1M^2}$, define 
\begin{equation*}
    \pi(I): A_{g_1g_2, g_2}(M^1)\otimes_{A_{g_2}(V)} M^2(0) \longmapsto M^3(0)
\end{equation*}
as 
\begin{equation*}
    \pi(I)(w_1\otimes w_2)=o_I(w_1)w_2,
\end{equation*}
where $w_1\in A_{g_1g_2, g_2}(M^1)$, $w_2\in M^2(0)$. Then by Lemma \ref{l3.8}, $$\pi (I) \in \mbox{Hom}_{A_{g_1g_2}(V)}(A_{g_1g_2, g_2}(M^1)\otimes_{A_{g_2}(V)} M^2(0), M^3(0)).$$
Thus we obtain a linear map
\begin{equation*}
    \pi: \mathcal{V}^{M^3}_{M^1M^2} \longmapsto \mbox{Hom}_{A_{g_1g_2}(V)}(A_{g_1g_2, g_2}(M^1)\otimes_{A_{g_2}(V)} M^2(0), M^3(0)).
\end{equation*}
\begin{theorem}\label{t3.13}
If $M^3$ is irreducible, then $\pi$ is injective. In particular, 
\[
\mbox{dim}\mathcal{V}^{M^3}_{M^1M^2} \leq \mbox{dim} \mbox{Hom}_{A_{g_1g_2}(V)}(A_{g_1g_2, g_2}(M^1)\otimes_{A_{g_2}(V)} M^2(0), M^3(0)).
\]
\end{theorem}
We need the following lemma to prove Theorem \ref{t3.13}. It's a generalization of Proposition \ref{p2.10}, Lemma \ref{l3.8} and
Proposition 4.5.7 in \cite{LL}.
\begin{lemma}\label{l3.14}
Let $I$ be an intertwining operator in $\mathcal{V}^{M^3}_{M^1M^2}$ such that $I^{\circ}(w_1, z)w_2=\sum_{n\in\frac{1}{T}\mathbb{Z}}w_1(n)w_2z^{-n-1}$. For $u\in V^{(j_1, j_2)}, w_1\in M^1$, $w_2\in M^2$, $n\in \frac{1}{T}\mathbb{Z}$ and $p\in \mathbb{Z}$, 
\[
u_{p+\frac{[j_1+j_2]}{T}}w_1(n)w_2=\sum_{i\in \Lambda}x_i(n_i)w_2,
\]
where $\Lambda$ is a finite index set, $x_i\in M^1$, $n_i\in \frac{1}{T}\mathbb{Z}$ and wt$(x_i(n_i))$=wt$(u_{p+\frac{[j_1+j_2]}{T}}w_1(n))$. In particular, if
wt$(u_{p+\frac{[j_1+j_2]}{T}}w_1(n))=0$, then $x_i(n_i)=o_I(x_i)$.
\end{lemma}
\begin{proof}
By associativity, let $k$ be an integer such that
\begin{align*}
&(z_0+z_2)^{k+\frac{j_2}{T}}Y_{M^3}(u, z_0+z_2)I^{\circ}(w_1, z_2)w_2\\
=&(z_2+z_0)^{k+\frac{j_2}{T}}I^{\circ}(Y_{M^1}(u, z_0)w_1, z_2)w_2
\end{align*}
Now we have
\begin{align*}
    &u_{p+\frac{[j_1+j_2]}{T}}w_1(n)w_2\\
    =&\mbox{Res}_{z_1}\mbox{Res}_{z_2}z_1^{p+\frac{[j_1+j_2]}{T}}z_2^{n}Y_{M^3}(u, z_1)I^{\circ}(w_1, z_2)w_2\\
    =&\mbox{Res}_{z_0}\mbox{Res}_{z_1}\mbox{Res}_{z_2}z_0^{-1}\delta(\frac{z_1-z_2}{z_0})z_1^{p+\frac{[j_1+j_2]}{T}}z_2^{n}Y_{M^3}(u, z_1)I^{\circ}(w_1, z_2)w_2\\
    \intertext{since $z_1^{p+\frac{[j_1+j_2]}{T}}Y_{M^3}(u, z_1)\in$ End$(M^3)[[z_1^{-1}, z_1]]$,}    
    =&\mbox{Res}_{z_0}\mbox{Res}_{z_2}(z_0+z_2)^{p+\frac{[j_1+j_2]}{T}}z_2^n Y_{M^3}(u, z_0+z_2)I^{\circ}(w_1, z_2)w_2\\
    =&\mbox{Res}_{z_0}\mbox{Res}_{z_2}(z_0+z_2)^{p-k+\frac{[j_1+j_2]}{T}-\frac{j_2}{T}}(z_0+z_2)^{k+\frac{j_2}{T}}z_2^n Y_{M^3}(u, z_0+z_2)I^{\circ}(w_1, z_2)w_2\\
    \intertext{let $k'$ be an integer such that $z_2^{k'+n}I^{\circ}(w_1, z_2)w_2\in z_2^{-1+\frac{1}{T}}M^3[[z_2^{\frac{1}{T}}]]$,}
    =&\mbox{Res}_{z_0}\mbox{Res}_{z_2}\sum_{i\geq 0}^{k'-1}\binom{p-k+\frac{[j_1+j_2]}{T}-\frac{j_2}{T}}{i}z_0^{p-k+\frac{[j_1+j_2]}{T}-\frac{j_2}{T}-i}z_2^i\\
    &\times z_2^n(z_0+z_2)^{k+\frac{j_2}{T}} Y_{M^3}(u, z_0+z_2)I^{\circ}(w_1, z_2)w_2\\
    =&\mbox{Res}_{z_0}\mbox{Res}_{z_2}\sum_{i\geq 0}^{k'-1}\binom{p-k+\frac{[j_1+j_2]}{T}-\frac{j_2}{T}}{i}z_0^{p-k+\frac{[j_1+j_2]}{T}-\frac{j_2}{T}-i}z_2^{i+n}\\
    &\times (z_2+z_0)^{k+\frac{j_2}{T}}I^{\circ}(Y_{M^1}(u, z_0)w_1, z_2)w_2\\
    =&\mbox{Res}_{z_0}\mbox{Res}_{z_2}\sum_{i\geq 0}^{k'-1}\sum_{j\geq 0}\binom{p-k+\frac{[j_1+j_2]}{T}-\frac{j_2}{T}}{i}\binom{k+\frac{j_2}{T}}{j}z_0^{p-k+\frac{[j_1+j_2]}{T}-\frac{j_2}{T}-i+j}\\
    &\times z_2^{i+n+k+\frac{j_2}{T}-j}I^{\circ}(Y_{M^1}(u, z_0)w_1, z_2)w_2\\
    =&\mbox{Res}_{z_0}\mbox{Res}_{z_2}\sum_{i\geq 0}^{k'-1}\sum_{j\geq 0}\binom{p-k+\frac{[j_1+j_2]}{T}-\frac{j_2}{T}}{i}\binom{k+\frac{j_2}{T}}{j}z_0^{p-k+\frac{[j_1+j_2]}{T}-\frac{j_2}{T}-i+j}\\
    &\times z_2^{i+n+k+\frac{j_2}{T}-j}\sum_{s\in \frac{1}{T}\mathbb{Z}}\sum_{t\in \frac{j_1}{T}+\mathbb{Z}}z_0^{-t-1}z_2^{-s-1}(u_tw_1)(s)w_2\\
    =&\sum_{i\geq 0}^{k'-1}\sum_{j\geq 0}\binom{p-k+\frac{[j_1+j_2]}{T}-\frac{j_2}{T}}{i}\binom{k+\frac{j_2}{T}}{j}\\
    &\times (u_{p-k+\frac{[j_1+j_2]}{T}-\frac{j_2}{T}-i+j}w_1)(i+n+k+\frac{j_2}{T}-j)w_2
\end{align*}
The proof is completed.
\end{proof}
Now we prove Theorem \ref{t3.13}.
\begin{proof}
It suffices to show that $\pi(I)=0$ implies $I^{\circ}=0$. Suppose $\pi (I)=0$, then 
\[
\langle f, I^{\circ}(w_1, z)w_2\rangle=0,\qquad \mbox{for any}\ f\in (M^3(0))^*,\ w_1\in M^1,\ w_2\in M^2(0).
\]
Furthermore, by Proposition \ref{p2.16},
\begin{align*}
\langle u_nf, I^{\circ}(w_1, z)w_2\rangle=&\langle u(n)f, I^{\circ}(w_1, z)w_2\rangle\\
=&\langle f, \theta(u(n))I^{\circ}(w_1, z)w_2\rangle
\end{align*}
here $u\in V^{(j_1, j_2)}$, $0\leq j_1, j_2\leq T-1$, $n\in \frac{[j_1+j_2]}{T}+\mathbb{Z}$ and  $u(n) \in V[g_1g_2]$.
By Lemma \ref{l3.14},
\[\langle f, \theta(u(n))I^{\circ}(w_1, z)w_2\rangle\] can be written as a sum of 
\[
\langle f, o_I(x_i)w_2\rangle,
\]
where $x_i\in M^1$ and $i\in \Lambda$, a finite index set. Thus we have 
\[
\langle f, \theta(u(n))I^{\circ}(w_1, z)w_2\rangle=0.
\]
Since $M^3$ is irreducible, $(M^3)'$ is generated by $(M^3(0))^*$. Therefore 
\[
\langle f, I^{\circ}(w_1, z)w_2\rangle=0,\qquad \mbox{for any}\ f\in (M^3)',\ w_1\in M^1,\ w_2\in M^2(0).
\]
Hence, by Lemma \ref{l3.14} again, we can prove the following stronger result:
\[
\langle f, u_{p+\frac{[j_1+j_2]}{T}}I^{\circ}(w_1, z)w_2\rangle=0
\]
for any $f\in (M^3)'$, $w_1\in M^1$, $w_2\in M^2(0)$, $u\in V^{(j_1, j_2)}$, $0\leq j_1, j_2\leq T-1$ and $p\in \mathbb{Z}$.\par
By assumption, $M^2$ is generated by $M^2(0)$. It remains to show 
\[
\langle f, I^{\circ}(w_1, z)u_{p+\frac{j_2}{T}}w_2\rangle=0
\]
for any $f\in (M^3)'$, $w_1\in M^1$, $w_2\in M^2(0)$, $u\in V^{(j_1, j_2)}$, $0\leq j_1, j_2\leq T-1$ and $p\in \mathbb{Z}$. Applying Res$_{z_0}$ to the generalized Jacobi identity for $I^{\circ}$, we get 
\begin{align*}
    &z_1^{\frac{j_2}{T}}(z_1-z_2)^{\frac{j_1}{T}}Y_{M^3}(u, z_1)I^{\circ}(w_1, z_2)w_2\\
    -&e^{\frac{j_1}{T}\pi i}z_1^{\frac{j_2}{T}}(z_2-z_1)^{\frac{j_1}{T}}I^{\circ}(w_1, z_2)Y_{M^2}(u, z_1)w_2\\
    =&\mbox{Res}_{z_0}z_0^{\frac{j_1}{T}}z_1^{-1}\delta(\frac{z_2+z_0}{z_1})(z_2+z_0)^{\frac{j_2}{T}}I^{\circ}(Y_{M^1}(u, z_0)w_1, z_2)w_2
\end{align*}
For any $p\in \mathbb{Z}$, applying Res$_{z_1}z_1^p$ to the above identity, we get  
\begin{align*}
    &\sum_{i\geq 0}\binom{\frac{j_1}{T}}{i}(-1)^iz_2^{\frac{j_1}{T}-i}I^{\circ}(w_1, z_2)u_{p+\frac{j_2}{T}+i}w_2\\
    =&e^{-\frac{j_1}{T}\pi i}\sum_{i\geq 0}\binom{\frac{j_1}{T}}{i}(-z_2)^iu_{p+\frac{j_1+j_2}{T}-i}I^{\circ}(w_1, z_2)w_2\\
    &-e^{-\frac{j_1}{T}\pi i}\sum_{i\geq 0}\binom{p+\frac{j_2}{T}}{i}z_2^{p+\frac{j_2}{T}-i}I^{\circ}(u_{i+\frac{j_1}{T}}w_1, z_2)w_2
\end{align*}
Thus
\[
\langle f, \sum_{i\geq 0}\binom{\frac{j_1}{T}}{i}(-1)^iz_2^{\frac{j_1}{T}-i}I^{\circ}(w_1, z_2)u_{p+\frac{j_2}{T}+i}w_2\rangle=0,
\]
for any $f\in (M^3)'$, $w_1\in M^1$, $w_2\in M^2(0)$, $u\in V^{(j_1, j_2)}$, $0\leq j_1, j_2\leq T-1$ and $p\in \mathbb{Z}$. Now let $k$ be the integer such that $u_{k+\frac{j_2}{T}+i}w_2=0$ for $i>0$, and choose $p$ to be $k$, then we get 
\[
\langle f, I^{\circ}(w_1, z_2)u_{k+\frac{j_2}{T}}w_2\rangle=0.
\]
Choose $p$ to be $k-1$, we get 
\[
\langle f, I^{\circ}(w_1, z_2)u_{k-1+\frac{j_2}{T}}w_2\rangle=0.
\]
Now we can inductively prove 
\[
\langle f, I^{\circ}(w_1, z_2)u_{p+\frac{j_2}{T}}w_2\rangle=0
\]
for any $p\in \mathbb{Z}$. The proof is completed.
\end{proof}

\begin{remark}
There is no reason for us not to conjecture that $\pi$ is surjective, hence $\mbox{dim}\mathcal{V}^{M^3}_{M^1M^2} = \mbox{dim} \mbox{Hom}_{A_{g_1g_2}(V)}(A_{g_1g_2, g_2}(M^1)\otimes_{A_{g_2}(V)} M^2(0), M^3(0))$. In \cite{L2} and \cite{L3}, two approaches were given for untwisted case when $g_1=g_2=g_3=1$. For twisted case, the challenge is that we can't prove $O'_{g_1g_2, g_2}(M^1)=O_{g_1g_2, g_2}(M^1)$. For  the approach in \cite{L3}, one of the key step is to identify $A(M^1)^*$ as $\Omega(\mathcal{D}_{P(-1)}(M^1))$, here $\mathcal{D}_{P(-1)}(M^1))$ is a $V\otimes V$-module in $(M^1)'^*$. But without condition $O'_{g_1g_2, g_2}(M^1)=O_{g_1g_2, g_2}(M^1)$, it is hard to make this identification in twisted case.     
\end{remark}
\section{Acknowledgement}
I am grateful to Professor Qifen Jiang at Shanghai Jiaotong University for valuable discussions.


\begin{thebibliography}{ABCD}






\bibitem[B]{B}
R. E. Borcherds, Vertex algebras, Kac-Moody algebras, and the Monster,
\emph{ Proc. Natl. Acad. Sci. USA} {\bf 83} (1986), 3068-3071.






\bibitem[DJ]{DJ}
C. Dong, C. Jiang, Bimodules associated to vertex operator algebra, \emph{Math.Z} \textbf{289} (2008), 799-826

\bibitem[DLM1]{DLM96}C. Dong, H. Li, G. Mason, Simple Currents and
Extensions of Vertex Operator Algebras. \emph{Comm. Math. Phys.} \textbf{180
}(1996), 671--707.

\bibitem[DLM2]{DLM97} C. Dong, H. Li, G. Mason, Regularity of rational
vertex operator algebras. \emph{Adv. Math.} \textbf{132} (1997), 148--166.

\bibitem[DLM3]{DLM98} C. Dong, H. Li, G. Mason, Twisted representations
of vertex operator algebras. \emph{Math. Ann. }\textbf{310 }(1998),
571--600.


\bibitem[DLM4]{DLM4}
C. Dong, H. Li and G. Mason, Vertex operator algebras and associative algebras, \emph{J. Alg.} \textbf{206} (1998), 67-96.

\bibitem[DR]{DR}
 C. Dong, L. Ren, Representations of vertex operator algebras and bimodules. \emph{J. Algebra} \textbf{384} (2013), 212–226.











\bibitem[FHL]{FHL} I. B. Frenkel, Y. Huang, J. Lepowsky, On axiomatic
approaches to vertex operator algebras and modules. \emph{Memoirs
American Math. Soc. }\textbf{104}, 1993.

\bibitem[FLM]{FLM} I. B. Frenkel, J. Lepowsky, A. Meurman, Vertex
operator algebras and the monster. \emph{Pure and Applied Math., }vol.
134, Academic Press, Massachusetts, 1988.

\bibitem[FZ]{FZ} I. Frenkel and Y. Zhu, Vertex operator algebras associated to
representations of affine and
Virasoro algebras, {\it Duke Math. J.} {\bf 66} (1992), 123-168.


\bibitem[JJ]{JJ}
Q. Jiang, X. Jiao, Bimodule and twisted representation of vertex operator algebras. \emph{Sci. China Math.} \textbf{59} (2016), no. 2, 397–410.

\bibitem[L1]{L1}
H. Li, Representation theory and tensor product theory for vertex  operator
algebras, Ph.D. thesis, Rutgers University, 1994.

\bibitem[L2]{L2}
H. Li, Determining fusion rules by $A(V)$-modules and bimodules, \emph{J. Alg.} \textbf{212} (1999), 515-556.

\bibitem[L3]{L3}
H. Li, The regular representations, Zhu’s $A(V)$-theory, and induced modules, \emph{J. Alg.} \textbf{238} (2001), 159-193.

\bibitem[L4]{L4}
H. Li, Regular representations of vertex operator algebras, \emph{Commun. Contemporary Math.} \textbf{4} (2002), 639-683.
\bibitem[LL]{LL}
J. Lepowsky, H. Li, Introduction to vertex operator algebras and their representations. In: \emph{Progress in Math}, Vol. 227. Boston: Birkhäuser, 2004

\bibitem[LS]{LS}
H. Li, J. Sun, Twisted regular representations of vertex operator
algebras, math-qa/2206.03455






\bibitem[MS]{MS} G. Moore and N. Seiberg, Classical and quantum
conformal field theory, {\em Comm. Math. Phys.} {\bf 123} (1989), 177-254.

\bibitem[MT]{MT}
M. Miyamoto, K. Tanabe, Uniform product of Ag,n(V) for an orbifold model V and G-twisted Zhu algebra, \emph{J. Alg.} \textbf{274} (2001): 80-96.

\bibitem[X]{X}X. Xu, Intertwining operators for twisted modules of a colored vertex operator
superalgebra, \emph{J. Alg.} \textbf{175} (1995), 241-273.
\bibitem[Z]{Z} Y. Zhu, Modular invariance of characters of vertex
operator algebras. \emph{J. Amer. Math. Soc.} \textbf{9} (1996), 237--302.
\end{thebibliography}
\end{document}